\documentclass{amsart}
\usepackage{epsfig}
\newtheorem{theorem}{Theorem}[section]
\newtheorem{prop}[theorem]{Proposition}
\newtheorem{lemma}[theorem]{Lemma}
\newtheorem{corollary}[theorem]{Corollary}
\newcommand\beq{\begin{equation}}
\newcommand\eeq{\end{equation}}
\newcommand\bce{\begin{center}}
\newcommand\ece{\end{center}}
\newcommand\bea{\begin{eqnarray}}
\newcommand\eea{\end{eqnarray}}
\newcommand\ben{\begin{enumerate}}
\newcommand\een{\end{enumerate}}

\newcommand\us{\underset}
\newcommand\wt{\widetilde}

\newcommand\nn{\nonumber}
\newcommand\bs{\bigskip}
\newcommand\ms{\medskip}
\newcommand\wh{\widehat}
\newcommand\brr{\begin{array}}
\newcommand\err{\end{array}}
\newcommand\bt{\begin{tabular}}
\newcommand\et{\end{tabular}}
\renewcommand\S{{\mathcal S}}
\newcommand\I{{\mathcal I}}
\newcommand\D{{\mathcal D}}
\newcommand\brs{\psi}
\newcommand\kra{\varphi}
\newcommand\ct{\mathrm{ct}}

\newcommand\rt{\mathrm{rt}}

\newcommand\fp{\mathrm{fp}}
\newcommand\exc{\mathrm{exc}}
\newcommand\des{\mathrm{des}}

\newcommand\C{{\mathbf C}}
\newcommand\F{{\mathbf F}}
\renewcommand\u{{\mathbf u}}
\renewcommand\d{{\mathbf d}}
\newcommand\simf{\sim}
\newcommand\simfe{\approx}
\newcommand\dep{\mathrm{depth}}
\newcommand\tdzero{\mathrm{td}_0}
\newcommand\tdneg{\mathrm{td}_{<0}}
\newcommand\at{\mathrm{all tun}}
\newcommand\Seq{\mathrm{Seq}}

\newcommand\Pyr{\mathcal Pyr}
\newcommand\E{{\mathcal E}}
\newcommand\Ds{{\mathcal Ds}}

\newenvironment{abstrac}{%
         \small
        \begin{center}%
          {\bfseries {Abstract}\vspace{-.5em}}%
                 \end{center}%
        \quotation}

\title{Multiple pattern avoidance with respect to
fixed points and excedances}

\begin{document}
\maketitle

\begin{center}
{\large Sergi Elizalde} \\[4pt]
Department of Mathematics \\
MIT, Cambridge, MA 02139\\
sergi@math.mit.edu
\end{center}

\bs

\begin{abstrac}
\ms We study the distribution of the statistics `number of fixed
points' and `number of excedances' in permutations avoiding
subsets of patterns of length~3. We solve all the cases of
simultaneous avoidance of more than one pattern, giving generating
functions enumerating these two statistics. Some cases are
generalized to patterns of arbitrary length. For avoidance of one
single pattern we give partial results. We also describe the
distribution of these statistics in involutions avoiding any
subset of patterns of length~3.

The main technique is to use bijections between pattern-avoiding
permutations and certain kinds of Dyck paths, in such a way that
the statistics in permutations that we study correspond to
statistics on Dyck paths that are easy to enumerate.
\end{abstrac}

\ms \bce MR Subject Classifications: 05A15, 05A05 \ece

\section{Introduction}

The problem of enumerating pattern-avoiding permutations, also
known as restricted permutations, has generated a lot of research
over the last few decades. One of the most referenced papers on
this topic is \cite{SS85}, which contains a systematic enumeration
of permutations avoiding any subset of patterns of length~3.

However, the study of statistics in pattern-avoiding permutations
started developing very recently, and the interest in this topic
is currently growing. Two of the most studied permutation
statistics have been the number of fixed points and the number of
excedances of a restricted permutation. For example, in
\cite{RSZ02,Eli02,EliPak} it is shown the surprising fact that the
joint distribution of these two statistics is the same in 132- and
in 321-avoiding permutations. In \cite{DRS}, involutions avoiding
any single pattern of length~3 are studied with respect to the
number of fixed points. Another paper \cite{ManRob02} deals with
the enumeration of permutations with a given number of fixed
points avoiding simultaneously two or more patterns of length~3.
Finally, \cite{GM02} considers additional restrictions on
132-avoiding involutions.

\ms

Given two permutations $\pi=\pi_1\pi_2\cdots\pi_n\in\S_n$ and
$\sigma=\sigma_1\sigma_2\cdots\sigma_m\in\S_m$, with $m\le n$, we
say that $\pi$ \emph{contains} $\sigma$ if there exist indices
$i_1<i_2<\ldots<i_m$ such that $\pi_{i_1}\pi_{i_2}\cdots\pi_{i_m}$
is in the same relative order as $\sigma_1\sigma_2\cdots\sigma_m$.
If $\pi$ does not contain $\sigma$, we say that $\pi$
\emph{avoids} $\sigma$, or that it is \emph{$\sigma$-avoiding}.
Denote by $\S_n(\sigma)$ the set of $\sigma$-avoiding permutations
in $\S_n$. More generally, if $\Sigma\subseteq\bigcup_{k\ge
1}\S_k$ is any set of permutations, denote by $\S_n(\Sigma)$ the
set of permutations in $\S_n$ that avoid simultaneously all the
patterns in $\Sigma$ (also called $\Sigma$-avoiding permutations).

We say that $i$ is a \emph{fixed point} of a permutation $\pi$ if
$\pi_i=i$, and that $i$ is an \emph{excedance} of $\pi$ if
$\pi_i>i$. Denote by $\fp(\pi)$ and $\exc(\pi)$ the number of
fixed points and the number of excedances of $\pi$ respectively.
We are interested in enumerating pattern-avoiding permutations
with respect to these two statistics. Most of the results will be
expressed in terms of multivariate generating functions.

\ms

In Section~\ref{sec:bijsta} we introduce the definitions, the
preliminaries, and the basic tools that we use in the paper.
Section~\ref{sec:single} is devoted to the study of permutations
with a single restriction with respect to the statistics $\fp$ and
$\exc$ (except in the case of the pattern 123, for which we can
only give partial results regarding $\fp$). In
Section~\ref{sec:double} we solve completely the case of
permutations avoiding simultaneously any two patterns of length~3,
giving generating functions counting the number of fixed points
and the number of excedances. For some particular cases we can
generalize the results, allowing one pattern of the pair to have
arbitrary length. In Section~\ref{sec:triple} we give the
analogous generating functions for permutations avoiding
simultaneously any three patterns of length~3 or more. Finally,
Section~\ref{sec:invol} is concerned with the study of the
distribution of these statistics in involutions avoiding any
subset of patterns of length~3. We end with a few final remarks
about the results of the paper and possible extensions.

\section{Bijections and statistics}\label{sec:bijsta}

We will represent a permutation $\pi\in\S_n$ by an $n\times n$
array with a cross in each one of the squares $(i,\pi_i)$. Fixed
points of $\pi$ correspond to crosses on the main diagonal (from
top-left to bottom-right corner) of the array, and excedances of
$\pi$ are represented by crosses to the right of this diagonal.

\subsection{Trivial bijections in $\S_n$}

Given a permutation $\pi=\pi_1\pi_2\cdots\pi_n$, define its
complementation $\bar\pi=(n+1-\pi_1)(n+1-\pi_2)\cdots(n+1-\pi_n)$.
The array of $\bar\pi$ is obtained from the array of $\pi$ by a
flip along a vertical axis, so fixed points (resp. excedances) of
$\pi$ correspond to crosses on (resp. to the left of) the
secondary diagonal (from top-right to bottom-left corner) of the
array of $\bar\pi$. Similarly, define $\wh\pi$ to be the
permutation whose array is the one obtained from that of $\pi$ by
reflection along the secondary diagonal. Note that reflecting the
array of $\pi$ along the main diagonal, we get the array of its
inverse $\pi^{-1}$. For any set of permutations $\Sigma$, let
$\bar\Sigma$ be the set obtained by reversing each of the elements
of $\Sigma$. Define $\wh\Sigma$ and $\Sigma^{-1}$ analogously. The
following trivial lemma will be used later on.

\begin{lemma}\label{lemma:refl} Let $\Sigma\in\bigcup_{k\ge 1}\S_k$ be any set of permutations, and let
$\pi\in\S_n$. We have that
\ben
\item $\pi\in\S_n(\Sigma) \iff \bar\pi\in\S_n(\bar\Sigma) \iff \wh\pi\in\S_n(\wh\Sigma) \iff \pi^{-1}\in\S_n(\Sigma^{-1})$,
\item $\fp(\wh\pi)=\fp(\pi)$, $\exc(\wh\pi)=\exc(\pi)$,
\item $\fp(\pi^{-1})=\fp(\pi)$,
$\exc(\pi^{-1})=n-\fp(\pi)-\exc(\pi)$. \een
\end{lemma}

We are interested in the distribution of the statistics $\fp$ and
$\exc$ among the permutations avoiding a certain set $\Sigma$ of
patterns. Given any such set $\Sigma$, we define the generating
function (GF for short) $F_\Sigma$ as
$$F_\Sigma(x,q,z):=\sum_{n\ge0}\sum_{\pi\in\S_n(\Sigma)}x^{\fp(\pi)}q^{\exc(\pi)}z^n.$$
If $\Sigma=\{\sigma\}$, we will write $F_\sigma$ instead of
$F_{\{\sigma\}}$. The following lemma restates the previous one in
terms of generating functions.

\begin{lemma}\label{lemma:sim} Let $\Sigma$ be any set of
permutations. We have \ben
\item $F_{\wh\Sigma}(x,q,z)=F_{\Sigma}(x,q,z)$,
\item
$F_{\Sigma^{-1}}(x,q,z)=F_{\Sigma}(\dfrac{x}{q},\dfrac{1}{q},qz)$.
\een
\end{lemma}

\begin{proof} To prove (1), consider the bijection between
$\S_n(\Sigma)$ and $\S_n(\wh\Sigma)$ that maps $\pi$ to $\wh\pi$.
The equation follows from parts (1) and (2) of
Lemma~\ref{lemma:refl}.

Equation (2) follows similarly from parts (1) and (3) of the
previous lemma, noticing that
\bea\nn\sum_{n\ge0}\sum_{\pi\in\S_n(\Sigma^{-1})}x^{\fp(\pi)}q^{\exc(\pi)}z^n=
\sum_{n\ge0}\sum_{\pi\in\S_n(\Sigma)}x^{\fp(\pi^{-1})}q^{\exc(\pi^{-1})}z^n\hspace{3cm}\\
\nn
=\sum_{n\ge0}\sum_{\pi\in\S_n(\Sigma)}x^{\fp(\pi)}q^{n-\fp(\pi)-\exc(\pi)}z^n
=\sum_{n\ge0}\sum_{\pi\in\S_n(\Sigma)}\left(\frac{x}{q}\right)^{\fp(\pi)}
\left(\frac{1}{q}\right)^{\exc(\pi)}(qz)^n.\eea
\end{proof}

If for two sets of patterns $\Sigma_1$ and $\Sigma_2$ we have that
$F_{\Sigma_1}(x,q,z)=F_{\Sigma_2}(x,q,z)$ (i.e., the joint
distribution of fixed points and excedances is the same in
$\Sigma_1$-avoiding as in $\Sigma_2$-avoiding permutations), as in
part (1) of Lemma~\ref{lemma:sim}, we will write
$\Sigma_1\simfe\Sigma_2$. If we are in the situation of part (2),
that is,
$F_{\Sigma_1}(x,q,z)=F_{\Sigma_2}(\frac{x}{q},\frac{1}{q},qz)$, we
will write $\Sigma_1\simf\Sigma_2$.

\subsection{Bijection to Dyck paths}\label{sec:bij}

Recall that a \emph{Dyck path} of length $2n$ is a lattice path in
$\mathbb{Z}^2$ between $(0,0)$ and $(2n,0)$ consisting of up-steps
$(1,1)$ and down-steps $(1,-1)$ which never goes below the
$x$-axis. We shall denote by $\D_n$ the set of Dyck paths of
length $2n$, and by $\D=\bigcup_{n\geq0}\D_n$ the class of all
Dyck paths. One of the main ingredients we will need throughout
the paper is a bijection between $\S_n(132)$ and $\D_n$, which we define
in the next paragraph.

\begin{figure}[hbt]
\epsfig{file=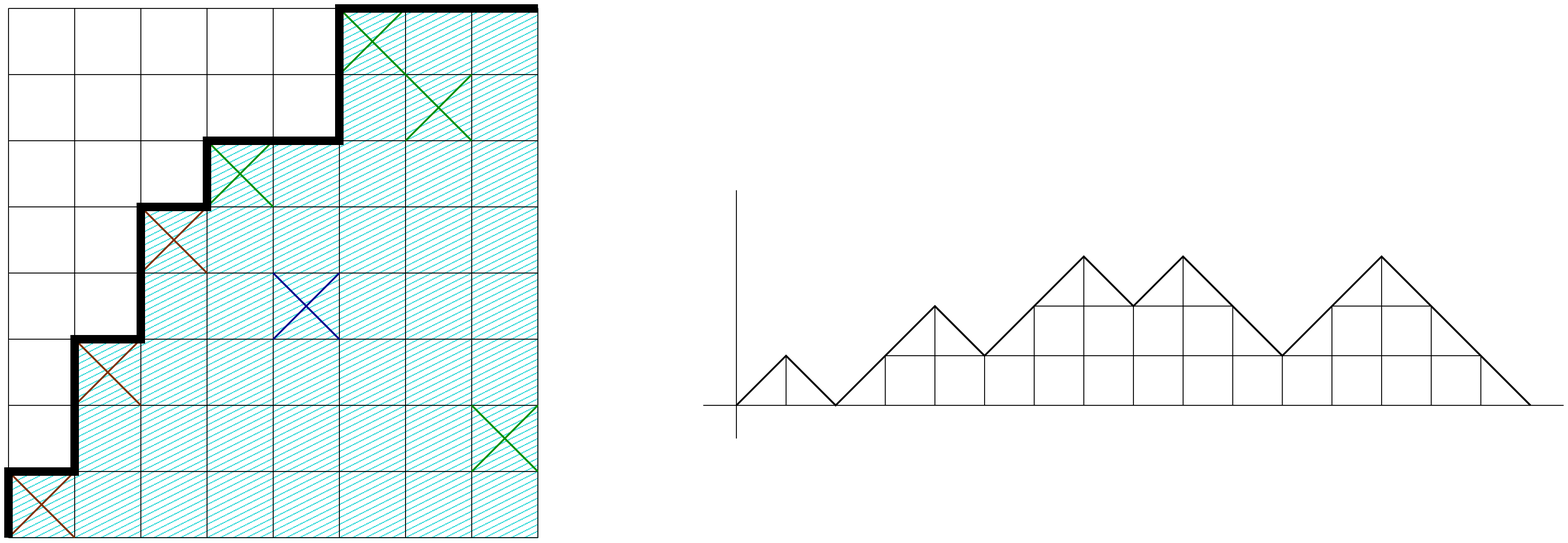,height=3.7cm} \caption{\label{fig:bij2}
The bijection $\kra$.}
\end{figure}

Any permutation $\pi\in\S_n$ can be represented as an $n\times n$ array
with crosses in positions $(i,\pi_i)$. From this array of crosses, we obtain the
\emph{diagram} of $\pi$ as follows. For each cross, shade the cell
containing it and the squares that are due south and due east of
it. The diagram is defined as the region that is left unshaded.
It is shown in~\cite{Rei02} that this gives a bijection
between $\S_n(132)$ and Young diagrams that fit in the shape
$(n-1,n-2,\ldots,1)$.
Consider now the path determined by the border of the diagram of $\pi$,
that is, the path with \emph{up} and \emph{right} steps that goes
from the lower-left corner to the upper-right corner of the array,
leaving all the crosses to the right, and staying always as close
to the diagonal connecting these two corners as possible. Define
$\kra(\pi)$ to be the Dyck path obtained from this path by
reading an up-step every time it goes up and a down-step every
time it goes right. Since the path in the array does not go below
the diagonal, $\kra(\pi)$ does not go below the $x$-axis.
Figure~\ref{fig:bij2} shows an example when $\pi=67435281$.

The bijection~$\kra$ is essentially the same bijection
between $\S_n(132)$ and $\D_n$ given by Krattenthaler~\cite{Kra01}
(see also~\cite{Ful01}), up to reflection of the path from a
vertical line.

\ms

Next we define the inverse map
$\kra^{-1}:\D_n\longrightarrow\S_n(132)$. Given a Dyck path
$D\in\D_n$, the first step needed to reverse the above procedure
is to transform $D$ into a path $U$ from the lower-left corner to
the upper-right corner of an $n\times n$ array, not going below
the diagonal connecting these two corners. Then, the squares to
the left of this path form a diagram, and we can shade all the
remaining squares. From this diagram, the permutation
$\pi\in\S_n(132)$ can be recovered as follows: row by row, put a
cross in the leftmost shaded square such that there is exactly one
cross in each column.

\subsection{Statistics on Dyck paths}

It is well-known that $|\D_n|=\C_n=\frac{1}{n+1}{2n \choose n}$,
the $n$-th Catalan number. If $D\in\D_n$, we will write $|D|=n$ to
indicate the {\em semilength} of $D$. The GF that enumerates Dyck
paths according to their semilength is
$\sum_{D\in\D}{z^{|D|}}=\sum_{n\geq0}{\C_n z^n}
=\frac{1-\sqrt{1-4z}}{2z}$, which we denote by $\C(z)$. Sometimes
it will be convenient to encode each up-step by a letter $\u$ and
each down-step by $\d$, obtaining an encoding of the Dyck path as
a \emph{Dyck word}.

A \emph{peak} of a Dyck path $D\in\D$ is an up-step followed by a
down-step (i.e., an occurrence of $\u\d$ in the associated Dyck
word). A \emph{hill} is a peak at height~1, where the height is
the $y$-coordinate of the top of the peak.
A \emph{valley} of a Dyck path $D\in\D$ is a down-step followed by
an up-step (i.e., an occurrence of $\d\u$ in the associated Dyck
word). The \emph{height} of $D$ is the $y$-coordinate of the
highest point of the path. Denote by $\D^{\le k}$ the set of Dyck
paths of height at most $k$.
Define a \emph{pyramid} to be a Dyck path that has only one peak,
that is, a path of the form $\u^k\d^k$ with $k\ge1$ (here the
exponent indicates the number of times the letter is repeated).

As defined in~\cite{Eli02}, a \emph{tunnel} of a Dyck path $D$ is
a horizontal segment between two lattice points of $D$ that
intersects $D$ only in these two points, and stays always below
$D$. Tunnels are in obvious one-to-one correspondence with
decompositions of the Dyck word $D=A\u B\d C$, where $B\in\D$ (no
restrictions on $A$ and $C$). In the decomposition, the tunnel is
the segment that goes from the beginning of the $\u $ to the end
of the $\d$. If $D\in\D_n$, then $D$ has exactly $n$ tunnels,
since such a decomposition can be given for each up-step $\u $ of
$D$. The \emph{length} of a tunnel is just its length as a
segment, and the \emph{height} is its $y$-coordinate. It will be
useful to define the \emph{depth} of a tunnel $T$ as
$\dep(T):=\frac{1}{2}\mathrm{length}(T)-\mathrm{height}(T)-1$.

A tunnel of $D\in\D_n$ is called a \emph{centered tunnel} if the
$x$-coordinate of its midpoint is $n$, that is, the tunnel is
centered with respect to the vertical line through the middle of
$D$. In terms of the decomposition $D=A\u B\d C$, this is
equivalent to saying that $A$ and $C$ have the same length. Denote
by $\ct(D)$ the number of centered tunnels of $D$.

\begin{figure}[hbt]
\epsfig{file=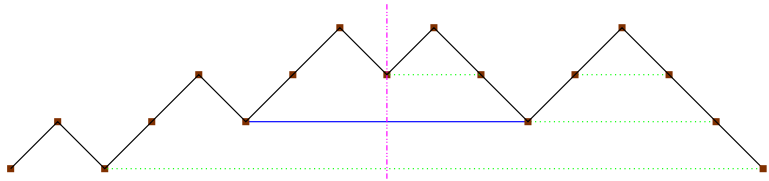,width=4in} \caption{\label{fig:ctrt} One
centered and four right tunnels.}
\end{figure}

A tunnel of $D\in\D_n$ is called a \emph{right tunnel} if the
$x$-coordinate of its midpoint is strictly greater than $n$, that
is, the midpoint of the tunnel is to the right of the vertical
line through the middle of $D$. Clearly, in terms of the
decomposition $D=A\u B\d C$, this is equivalent to saying that the
length of $A$ is strictly bigger than the length of $C$. Denote by
$\rt(D)$ the number of right tunnels of $D$. In
Figure~\ref{fig:ctrt}, there is one centered tunnel drawn with a
solid line, and four right tunnels drawn with dotted lines.

For a Dyck path $D\in\D_n$, denote by $D^\ast$ the path obtained
by reflection of $D$ from a vertical line $x=n$. We say that $D$
is {\em symmetric} if $D=D^\ast$. Denote by $\Ds\subset\D$ the
subclass of symmetric Dyck paths.

\subsection{Properties of $\kra$}

In order to enumerate fixed points and excedances in permutations,
we analyze what these statistics are mapped to by $\kra$.
Table~\ref{prokra} summarizes the correspondences of $\kra$ that
we will use.

\begin{table}[htb]
  \bt{c|c|c} {\it In the permutation $\pi$} & {\it In the array of $\pi$} & {\it In the Dyck path $\kra(\pi)$}
  \\ \hline
  fixed points of $\pi$ & crosses on the main diagonal & centered
  tunnels \\ \hline
  excedances of $\pi$ & crosses to the right of the main diagonal & right
  tunnels \\ \hline
  fixed points of $\bar\pi$ & crosses on the secondary diagonal & tunnels of depth 0
  \\ \hline
  excedances of $\bar\pi$ & crosses to the left of the secondary diagonal &
  tunnels of depth $<0$ \\
  \et\ms
  \caption{\label{prokra} Behavior of $\kra$ on fixed points and excedances.}
\end{table}
The correspondences between the first two columns have been
discussed above. Now we show how $\kra$ maps the statistics on the
second column to those on the third one. We repeat the reasoning
in \cite{Eli02,EliPak}.

For this purpose, instead of using $D=\kra(\pi)$, it
will be convenient to consider the path $U$ from the lower-left corner
to the upper-right corner of the array of $\pi$. We will talk
about tunnels of $U$ to refer to the corresponding tunnels of $D$
under this trivial transformation.

\begin{figure}[hbt]
\epsfig{file=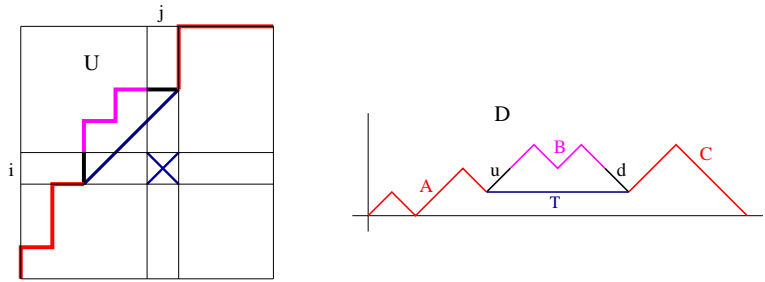,height=3.7cm} \caption{\label{fig:tuncross}
Each cross produces a tunnel.}
\end{figure}

Consider the arrangement of crosses of $\pi$ as defined earlier. We now show how to associate
a unique tunnel $T$ of $D$ to each cross $X$ of this array. Observe that given a cross
$X=(i,j)$, $U$ has a vertical step in row $i$ and a horizontal step in column $j$.
In $D$, these two steps in $U$ correspond to steps
$\u$ and $\d$ respectively, so they determine a decomposition $D=A\u B\d C$ (see Figure~\ref{fig:tuncross}),
and therefore a tunnel $T$ of $D$.

The distance between these two steps determines
the length of $T$, and the distance from these steps to the
secondary diagonal of the array determines the height of $T$. In
order for the corresponding cross to lie on the secondary
diagonal, the relation between these two quantities must be
$\frac{1}{2}\mathrm{length}(T)=\mathrm{height}(T)+1$, which is
equivalent to $\dep(T)=0$, by the definition of depth. The depth
of $T$ indicates how far from the secondary diagonal $X$ is. The
cross lies to the left of the secondary diagonal exactly when
$\dep(T)<0$. This justifies the last two rows of the table. The
first two correspondences are easier to see. Indeed, $T$ is
centered precisely when $X$ is on the main diagonal, and $T$ is a
right tunnel when $X$ lies to the right of the main diagonal.

\begin{figure}[hbt]
\epsfig{file=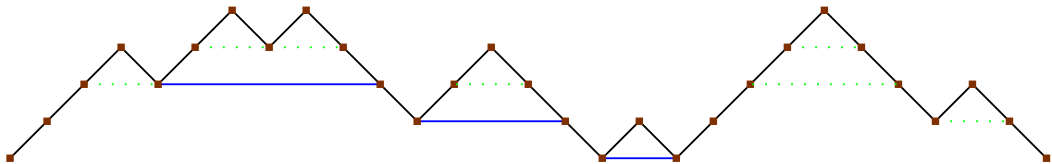,width=5.5in} \caption{\label{fig:depth}
Three tunnels of depth 0 and seven tunnels of depth $<0$.}
\end{figure}

We define two new statistics on Dyck paths. For $D\in\D$, let
$\tdzero(D)$ be the number of tunnels of depth 0 of $D$, and let
$\tdneg(D)$ be the number of tunnels of depth $<0$ of $D$. In
Figure~\ref{fig:depth}, there are three tunnels of depth 0 drawn
with a solid line, and seven tunnels of depth $<0$ drawn with
dotted lines. We have proved the following lemma.

\begin{lemma}\label{lemma:kra} Let $\pi\in\S_n(132)$, $\rho\in\S_n(312)$. We have
\ben\item $\fp(\pi)=\ct(\kra(\pi))$, \item
$\exc(\pi)=\rt(\kra(\pi))$,
\item $\fp(\rho)=\tdzero(\kra(\bar\rho))$, \item
$\exc(\rho)=\tdneg(\kra(\bar\rho))$. \een
\end{lemma}

\subsection{Combinatorial classes and generating functions.} Here
we direct the reader to \cite{FlSe98} for a detailed account on
combinatorial classes and the symbolic method. Let ${\mathcal A}$
be a class of unlabelled combinatorial objects and let $|\alpha|$
be the size of an object $\alpha \in {\mathcal A}$. If ${\mathcal
A}_n$ denotes the objects in ${\mathcal A}$ of size $n$ and $a_n =
|{\mathcal A}_n|$, then the {\em ordinary generating function} of
the class ${\mathcal A}$ is
$$ A(z) = \sum_{\alpha\in {\mathcal A}} z^{|\alpha|} = \sum_{n\ge0} a_n z^n.
 $$
In our context, the size of a Dyck path is simply its semilength.

There is a direct correspondence between set theoretic operations
(or ``constructions'') on combinatorial classes and algebraic
operations on GFs. Table~\ref{tab:symb} summarizes this
correspondence for the operations that are used in the paper.
There ``union'' means union of disjoint copies, ``product'' is the
usual cartesian product, and ``sequence'' forms an ordered
sequence in the usual sense.
\begin{table}[htb]
\begin{tabular}{l|l|l} \hline
\em Construction &  &  \em Operation on GFs \\
 \hline
 Union & ${\mathcal A}={\mathcal B}+{\mathcal C}$  &  $A(z)=B(z)+C(z)$ \\
 Product & ${\mathcal A}={\mathcal B}\times{\mathcal C}$  & $A(z)=B(z)C(z)$ \\
 Sequence & ${\mathcal A} = \Seq({\mathcal B})$ &  $A(z) = \frac{1}{1-B(z)}$ \\
 \hline
\end{tabular}
\ms \caption{\label{tab:symb} The basic combinatorial
constructions and their translation into GFs.}
\end{table}
Enumerations according to size and auxiliary parameters
$\chi_1,\chi_2,\ldots,\chi_r$ are described by multivariate GFs,
$$ A(u_1,u_2,\ldots,u_r,z)=\sum_{\alpha\in\mathcal A}
u_1^{\chi_1(\alpha)}u_2^{\chi_2(\alpha)}\cdots
u_r^{\chi_r(\alpha)}z^{|\alpha|}. $$

Throughout the paper the variable $z$ is reserved for marking the
length of a permutation and the semilength of a Dyck path, $x$ is
used for marking the number of fixed points of a permutation and
the number of centered tunnels or tunnels of depth 0 of a Dyck
path, and $q$ is the variable that marks the number of excedances
of a permutation and the number of right tunnels or tunnels of
depth $<0$ of a Dyck path.

\section{Single restrictions}\label{sec:single}

The most difficult case appears to be that of permutations
avoiding one single pattern. It is well known that for any
$\sigma\in\S_3$, $|\S_n(\sigma)|=\C_n$. By Lemma~\ref{lemma:sim},
we have that $132\simfe 213$, and that $231\simf 312$. These are
the only equivalences that follow from the trivial bijections.
Recently it was shown that, surprisingly, $321\simfe 132$. The
first proof of this result appears in \cite{Eli02} (see
\cite{EliPak} for a bijective proof). The weaker version $321\simf
132$ was proved earlier in \cite{RSZ02}.

So, we have the following equivalence classes of patterns of
length~3 with respect to fixed points and excedances:

\bce{\bf a) $123$\\
b) $132\simfe 213\simfe 321$ \\
c) $231 \quad\simf\quad$ c') $312$}\ece

\subsection{a) $123$}\label{s:123}

For this case we have not been able to find a satisfactory
expression for $F_{123}(x,q,z)$. We can nevertheless give
summation formulas for the number of permutations in $\S_n(123)$
with a given number of fixed points. The first trivial observation
is that if $\pi$ avoids 123, then it can have at most two fixed
points. If $\pi_i=i$, we say that $i$ is a \emph{big fixed point}
of $\pi$ if $i\ge\frac{n+1}{2}$, and that it is a \emph{small
fixed point} if $i<\frac{n+1}{2}$.

It is known that a permutation is 321-avoiding if and only if both
the subsequence determined by its excedances and the one
determined by the remaining elements are increasing. Using the
fact that $\pi$ avoids 123 if and only if $\bar\pi$ avoids 321, we
obtain a characterization of 123-avoiding permutations as those
with the following property: the elements $\pi_i$ such that
$\pi_i<n+1-i$ form a decreasing subsequence, and so do the
remaining elements. In particular, since no two fixed points can
be in the same decreasing subsequence, this implies that $\pi$ can
have at most one big fixed point and one small fixed point.

\begin{figure}[hbt]
\epsfig{file=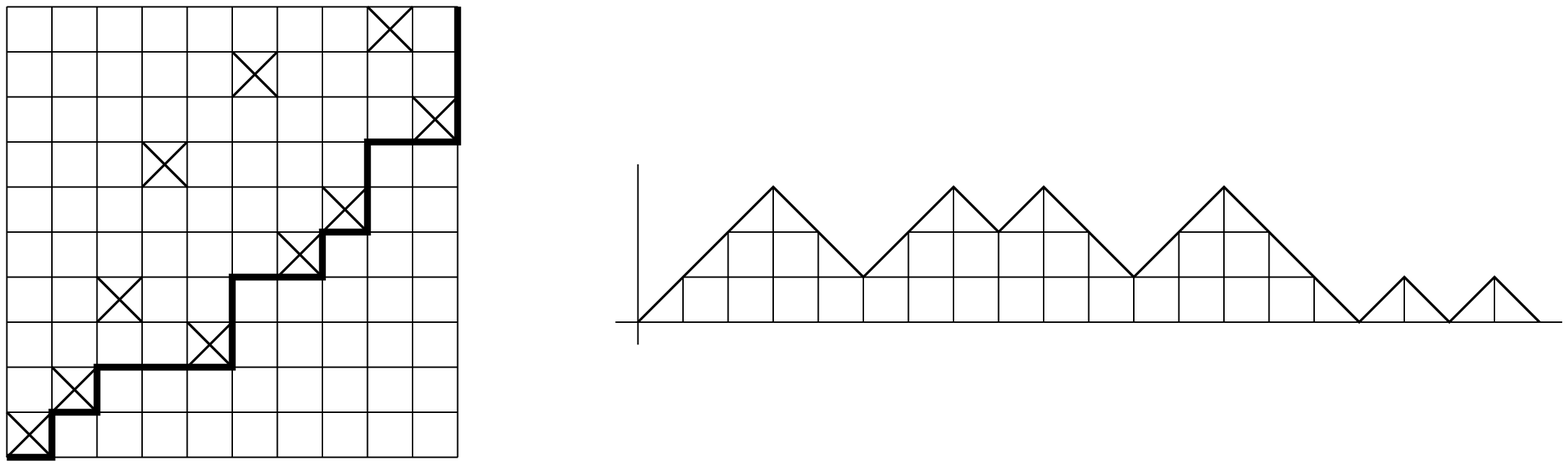,height=3.7cm} \caption{\label{fig:bij_brs}
The bijection $\brs$.}
\end{figure}

In this subsection we use a bijection between $\S_n(123)$ and
$\D_n$, which we denote by $\brs$. Up to composition with the
complementation operation, it is essentially the same as the
bijection from $\S_n(321)$ and $\D_n$ given in \cite{Eli02}. We
can define $\brs$ as follows. Let $\pi\in\S_n(123)$ be represented
as an $n\times n$ array with a cross on each square $(i,\pi_i)$.
Consider the path with \emph{down} and \emph{left} steps that goes
from the upper-right corner to the lower-left corner, leaving all
the crosses to the left, and staying always as close to the
diagonal connecting these two corners as possible. Then
$\brs(\pi)$ is the Dyck path obtained from this path by reading an
up-step every time the path goes down and a down-step every time
it goes left. Figure~\ref{fig:bij_brs} shows an example when
$\pi=(9,6,10,4,8,7,3,5,2,1)$.

Note that the peaks of the path are determined by the crosses of
elements $\pi_i$ such that $\pi_i\ge n+1-i$, which form a
decreasing subsequence. Now it is easy to determine how many
permutations have a big (resp. small) fixed point.

\begin{prop}\label{bigsmall} Let $n\ge1$. We have \ben \item $|\{\pi\in\S_n(123):\pi \mbox{ has a big fixed
point}\}|=\C_{n-1},$ \item $|\{\pi\in\S_n(123):\pi \mbox{ has a
small fixed
point}\}|=\left\{\bt{ll} $\C_{n-1}$ & if $n$ is even, \\
$\C_{n-1}-\C_{\frac{n-1}{2}}^2$ & if $n$ is odd. \et\right.$ \een
\end{prop}

\begin{proof}
(1) It is clear from the definition of $\brs$ that $\pi$ has a big
fixed point if and only if $\brs(\pi)$ has a peak in the middle.
Now, we can easily define a bijection from the subset of elements
of $\D_n$ with a peak in the middle and $\D_{n-1}$, by removing
the two middle steps $\u\d$.

(2) Clearly, $\pi\in\S_n(123)$ if and only if
$\wh\pi\in\S_n(123)$. This involution switches big and small fixed
points, except for the possible big fixed point in position
$\frac{n+1}{2}$, which remains unchanged. Applying now $\brs$, a
small fixed point of $\pi$ is transformed into a peak in the
middle of $\brs(\wh\pi)$ of height at least two (indeed, a hill
would correspond to the big fixed point $\frac{n+1}{2}$). Knowing
that the number of paths in $\D_n$ with a peak in the middle is
$\C_{n-1}$, we just have to subtract those where this middle peak
has height~1. If $n$ is even, paths in $\D_n$ cannot have a hill
in the middle. If $n$ is odd, such paths have the form $A\u\d B$,
where $A,B\in\D_{\frac{n-1}{2}}$, so the formula follows.
\end{proof}

For $k\ge0$, let $s_n^k(123):=|\{\pi\in\S_n(123):\fp(\pi)=k\}|$.
We have mentioned that $s_n^k(123)=0$ for $k\ge3$. The following
corollary reduces the problem of studying the distribution of
fixed points in $S_n(123)$ to that of determining $s_n^2(123)$.

\begin{corollary}\label{cor:123} Let $n\ge1$. We have \ben
\item $s_n^1(123)=\left\{\bt{ll} $2(\C_{n-1}-s_n^2(123))$ &
if $n$ even, \\ $2(\C_{n-1}-s_n^2(123))-\C_{\frac{n-1}{2}}^2$ & if
$n$ odd. \et\right. $
\item  $s_n^0(123)=\left\{\bt{ll} $\C_n-2\C_{n-1}+s_n^2(123)$ &
if $n$ even, \\ $\C_n-2\C_{n-1}+s_n^2(123)+\C_{\frac{n-1}{2}}^2$ &
if $n$ odd. \et\right. $ \een \end{corollary}

\begin{proof}
(1) By inclusion-exclusion, $s_n^1(123)=|\{\pi\in\S_n(123):\pi
\mbox{ has a big fixed point}\}|+|\{\pi\in\S_n(123):\pi \mbox{ has
a small fixed point}\}|-2s_n^2(123)$. Now we apply
Proposition~\ref{bigsmall}.

(2) Clearly, $s_n^0(123)=\C_n-s_n^1(123)-s_n^2(123)$.
\end{proof}

The next theorem, together with the previous corollary, gives a
formula for the distribution of fixed points in 123-avoiding
permutations.

\begin{theorem}\label{th:messy123}
\bea\nn s_n^2(123)&=&\sum_{i=1}^{n-1}\sum_{r,s=1}^{i}
\left[\left({2i-r-1\choose i-1}-{2i-r-1\choose
i}\right)\cdot\left({2i-s-1\choose i-1}-{2i-s-1\choose i}\right)\right.\\
\nn && \hspace*{2cm}\left.\cdot\us{n-h \ \mathrm{even}}
{\sum_{h=1}^n}\sum_{k=0}^{n-2i} f(k,r,h,n-2i+r)
f(n-2i-k,s,h,n-2i+s)\right],\eea where
$$f(k,r,h,\ell)=\left\{\bt{cl}
  $\dbinom{\frac{\ell+h-r}{2}-1}{k}\dbinom{\frac{\ell-h+r}{2}-1}{k-1}-
\dbinom{\frac{\ell-h-r}{2}-1}{k}\dbinom{\frac{\ell+h+r}{2}-1}{k-1}$
&
  if $k\ge 1$,\\
  $1$ & if $k=0$ and $\ell=h-r$,\\
  $0$ & otherwise,
\et\right.$$ with the convention ${a \choose b}:=0$ if $a<0$.
\end{theorem}

\begin{proof} Recall that $s_n^2(123)$ counts permutations with both
a big and a small fixed point. We have seen already that $\brs$
maps a big fixed point of the permutation into a peak in the
middle of the Dyck path. Now we look at how a small fixed point of
the permutation is transformed by $\brs$. We claim that
$\pi\in\S_n(123)$ has a small fixed point if and only if
$D=\brs(\pi)$ satisfies the following condition (which we call
condition C1): there exists $i$ such that the $i$-th and
$(i+1)$-st up-steps of $D$ are consecutive, the $i$-th and
$(i+1)$-st down-steps from the end are consecutive, and there are
exactly $n+1-2i$ peaks of $D$ between them. To see this, assume
that $i$ is a small fixed point of $\pi$ (see
Figure~\ref{fig:dpeaks}).
Then, the path from the upper-right to the lower-left corner of
the array of $\pi$, used to define $\brs(\pi)$, has two
consecutive vertical steps in rows $i$ and $i+1$, and two
consecutive horizontal steps in columns $i$ and $i+1$. Besides,
there are $n+1-2i$ crosses below and to the right of cross
$(i,i)$, each one of which produces a peak in the Dyck path
$\brs(\pi)$. Reciprocally, it can be checked that if $\brs(\pi)$
satisfies condition C1 then $\pi$ has a small fixed point.

\begin{figure}[hbt]
\epsfig{file=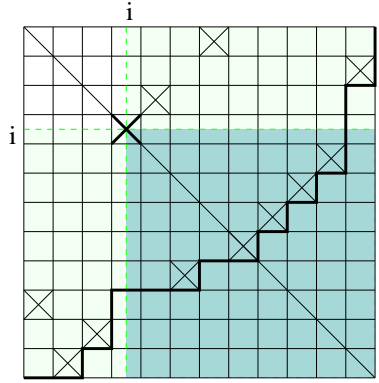,height=2in} \caption{\label{fig:dpeaks}
A small fixed point $i$ has $n+1-2i$ crosses below and to the right.}
\end{figure}

All we have to do is count how many paths $D\in\D_n$ with a peak
in the middle satisfy condition C1. For such a Dyck path $D$,
define the following parameters: let $i$ be the value such that
condition C1 holds, let $h$ be the height of $D$ in the middle,
$r$ the height at which the $i$-th up-step ends, and $s$ the
height of the point between the $i$-th and $(i+1)$-st down-steps
from the end. In the example of Figure~\ref{fig:param}, $n=12$,
$i=4$, $h=4$, $r=3$, and $s=1$.

\begin{figure}[hbt]
\epsfig{file=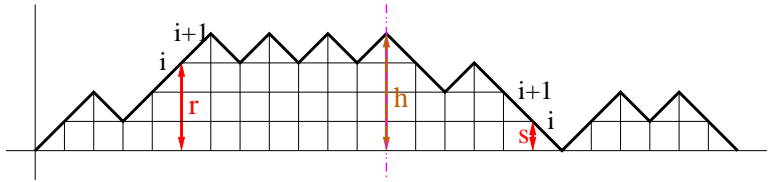,width=4in} \caption{\label{fig:param}
The parameters $i$, $h$, $r$ and $s$ in a Dyck path.}
\end{figure}

Fix $n$, $i$, $h$, $r$ and $s$. We will count the number of Dyck
paths $D$ with these given parameters. We can write $D=A\u\u
B_1B_2\d\d C$, where the distinguished $\u$'s are the $i$-th and
$(i+1)$-st up-steps, similarly with the $\d$'s, and the middle of
$D$ is between $B_1$ and $B_2$. Then $A$ is a path from $(0,0)$ to
$(2i-r-1,r-1)$ not going below $x=0$. It is easy to see that there
are ${2i-r-1\choose i-1}-{2i-r-1\choose i}$ such paths $A$. By
symmetry, there are ${2i-s-1\choose i-1}-{2i-s-1\choose i}$
possibilities for $C$.

Now we count the possibilities for $B_1$ and $B_2$. It can be
checked that $f(k,r,h,\ell)$ counts the number of paths from
$(0,r)$ to $(\ell,h)$ having exactly $k$ peaks, starting and
ending with an up-step, and never going below $x=0$. The fragment
$\u B_1$ is a path from $(2i-r,r)$ to $(n,h)$ not going below
$x=0$, and ending with an up-step (since $D$ has a peak in the
middle). If we fix $k$ as the number of peaks of this fragment,
then there are $f(k,r,h,n-2i+r)$ such paths $\u B_1$. Similarly,
there are $f(n-2i-k,s,h,n-2i-s)$ possibilities for $B_2\d$ with
$n-2i-k$ peaks.

Summing over all possible values of $k$, $h$, $r$, $s$ and $i$ we
obtain the expression in the theorem.
\end{proof}

\ms

Using Corollary~\ref{cor:123}, we can prove that
among the derangements (i.e., permutations without fixed points)
of length $n$, the number of 123-avoiding ones is at least the
number of 132-avoiding ones. This inequality was conjectured by
Mikl\'os B\'ona and Olivier Guibert.

\begin{theorem}[\cite{BG}]
For all $n\ge4$, $s_n^0(132)<s_n^0(123)$.
\end{theorem}
\begin{proof}
For $n\le12$ the result can be checked by exhaustive enumeration
of all derangements. Let us assume that $n\ge13$.

From part (1) of Corollary~\ref{cor:123}, we have that
$$s_n^1(123)+2s_n^2(123)=\left\{\bt{ll} $2\C_{n-1}$ & if $n$ even,
\\ $2\C_{n-1}-\C_{\frac{n-1}{2}}^2$ & if $n$ odd.
\et\right. $$ Thus, $$s_n^0(123)=\C_n-s_n^1(123)-s_n^2(123)\ge
\C_n-s_n^1(123)-2s_n^2(123)\ge\C_n-2\C_{n-1}.$$

It is known \cite{RSZ02} that $s_n^0(132)=\F_n$, the $n$-th {\em
Fine number}. Therefore, the theorem will be proved if we show
that \beq\label{ineqFC}\F_n<\C_n-2\C_{n-1}\eeq for $n\ge13$. Using
the identity
$\F_n=\frac{1}{2}\sum_{i=0}^{n-2}(\frac{-1}{2})^i\C_{n-i}$, we get
the inequality $\F_n< \frac{1}{2}\C_n - \frac{1}{4}\C_{n-1} +
\frac{1}{8}\C_{n-2}$, which reduces (\ref{ineqFC}) to showing that
$\C_n > \frac{7}{2}\C_{n-1} + \frac{1}{4}\C_{n-2}$. This
inequality certainly holds asymptotically, because $\C_n$ grows
like $\frac{1}{\sqrt{\pi}}n^{-\frac{3}{2}}4^n$ as $n$ tends to
infinity, and it is not hard to see that in fact it holds for all
$n\ge13$.
\end{proof}

\subsection{b) $132\simfe 213\simfe 321$}

This case has already been studied in \cite{Eli02}. The
corresponding GF, which appeared independently in \cite{Vel03} for
the case of the pattern 321, is the following.

\begin{theorem}[\cite{Eli02,Vel03}]\label{th:fpexc132}
$$F_{132}(x,q,z)=F_{213}(x,q,z)=F_{321}(x,q,z)=
\frac{2}{1+z(1+q-2x)+\sqrt{1-2z(1+q)+z^2(1-q)^2}}.$$
\end{theorem}

In fact, permutations avoiding these patterns have been enumerated
with respect to an additional statistic, the number of descents.
Recall that $i\le n-1$ is a \emph{descent} of $\pi\in\S_n$ if
$\pi_i>\pi_{i+1}$. Denote by $\des(\pi)$ the number of descents of
$\pi$. Theorem~\ref{th:fpexc132} can be generalized as follows.

\begin{theorem}[\cite{Eli02,EliDeu}]
The GF for 321-avoiding permutations with respect to fixed points,
excedances, and descents is
$$\sum_{n\geq0}\sum_{\pi\in\S_n(321)}x^{\fp(\pi)}q^{\exc(\pi)}p^{\des(\pi)}z^n
= \frac{2}{1+z(1+q-2x)+\sqrt{1-2z(1+q)+z^2((1+q)^2-4qp)}}.$$
 Similarly,

\bea\nn 1+\sum_{n\ge 1}
\sum_{\pi\in\S_n(132)}{x^{\fp(\pi)}q^{\exc(\pi)}}p^{\mathrm{des}(\pi)+1}z^n
= 1+\sum_{n\ge 1}
\sum_{\pi\in\S_n(213)}{x^{\fp(\pi)}q^{\exc(\pi)}}p^{\mathrm{des}(\pi)+1}z^n
\hspace{1cm} \\ \nn=\frac{2(1+xz(p-1))}{1+(1+q-2x)z-qz^2
(p-1)^2+\sqrt{f_1(q,z)}},\eea where $f_1(q,z)=
1-2(1+q)z+[(1-q)^2-2q(p-1)(p+3)]z^2-2q(1+q)(p-1)^2 z^3+q^2 (p-1)^4
z^4$.

\end{theorem}

\subsection{\bf c, c') $231\simf 312$}

Using the bijection \bce\bt{ccc} $\S_n(312)$ &
$\longleftrightarrow$ & $\D_n$ \\ $\pi$ & $\mapsto$ &
$\kra(\bar\pi)$, \et\ece Lemma~\ref{lemma:kra} implies that
$$F_{312}(x,q,z)=\sum_{D\in\D}x^{\tdzero(D)}q^{\tdneg(D)}z^{|D|}.$$
To enumerate tunnels of depth 0, we will separate them according
to their height. For every $h\ge 0$, a tunnel at height $h$ must
have length $2(h+1)$ in order to have depth 0. It is important to
notice that if a tunnel of depth 0 of $D$ corresponds to a
decomposition $D=A\u B\d C$, then $D$ has no tunnels of depth 0 in
the part given by $B$. In other words, the orthogonal projections
on the $x$-axis of all the tunnels of depth 0 of a given Dyck path
are disjoint. This observation allows us to give a continued
fraction expression for $F_{312}(x,1,z)$.

\begin{theorem} $F_{312}(x,1,z)$ is given by the following
continued fraction.
$$F_{312}(x,1,z)=\frac{1}{1-(x-1)z-\dfrac{z}{1-(x-1)z^2-\dfrac{z}{1-2(x-1)z^3-
\dfrac{z}{1-5(x-1)z^4-\dfrac{z}{\ddots}}}}},$$ where at each
level, the coefficient of $(x-1)z^{n+1}$ is the Catalan number
$\C_n$.
\end{theorem}

\begin{proof}
For every $h\ge0$, let $\tdzero^h(D)$ be the number of tunnels of
$D$ of height $h$ and length $2(h+1)$. Note that
$\tdzero(D)=\sum_{h\ge0}\tdzero^h(D)$. We will show now that for
every $h\ge1$, the generating function for Dyck paths where $x$
marks the statistic $\tdzero^0(\cdot)+\cdots+\tdzero^{h-1}(\cdot)$ is
given by the continued fraction of the theorem truncated at level
$h$, with the $(h+1)$-st level replaced with $\C(z)$.

A Dyck path $D$ can be written uniquely as a sequence of elevated
Dyck paths, that is, as $D=\u A_1\d\cdots\u A_r\d$, where each
$A_i\in\D$. In terms of the GF $\C(z)=\sum_{D\in\D}z^{|D|}$, this
translates into the equation $\C(z)=\frac{1}{1-z \C(z)}$. A tunnel
of height~0 and length~2 (i.e., a hill) appears in $D$ for each
empty $A_i$. Therefore, the GF enumerating hills is
\beq\sum_{D\in\D}x^{\tdzero^0(D)}z^{|D|}=\frac{1}{1-z[x-1+\C(z)]},\label{height0}\eeq
since an empty $A_i$ has to be accounted as $x$, not as 1.

Let us enumerate simultaneously hills (as above), and tunnels of
height~1 and length~4. The GF~(\ref{height0}) can be written as
$$\frac{1}{1-z\left[x-1+\dfrac{1}{1-z \C(z)}\right]}.$$
Combinatorially, this corresponds to expressing each $A_i$ as a
sequence $\u B_1\d\cdots\u B_s\d$, where each $B_j\in\D$. Notice
that since each $\u B_j\d$ starts at height~1, a tunnel of
height~1 and length~4 is created by each $B_j=\u\d$ in the
decomposition. Thus, if we want $x$ to mark also these tunnels,
such a $B_j$ has to be accounted as $xz$, not $z$. The
corresponding GF is
$$\sum_{D\in\D}x^{\tdzero^0(D)+\tdzero^1(D)}z^{|D|}=\frac{1}{1-z\left[x-1+\dfrac{1}{1-z
[(x-1)z+\C(z)]}\right]}.$$

Now it is clear how iterating this process indefinitely we obtain
the continued fraction of the theorem. From the GF where $x$ marks
$\tdzero^0(\cdot)+\cdots+\tdzero^{h-1}(\cdot)$, we can obtain the one
where $x$ marks $\tdzero^0(\cdot)+\cdots+\tdzero^h(\cdot)$ by replacing
the $\C(z)$ at the lowest level with
$$\frac{1}{1-z[(x-1)\C_hz^h+\C(z)]},$$ to account for tunnels of
height $h$ and length $2(h+1)$, which in the decomposition
correspond to elevated Dyck paths at height $h$.
\end{proof}

The same technique can be used to enumerate excedances in
312-avoiding permutations, which correspond to tunnels of depth
$<0$ in the Dyck path.

\begin{theorem}\label{th:312} Let $$\C_{<i}(z)=\sum_{j=0}^{i-1}\C_j z^j$$ be the
series for the Catalan numbers truncated at degree $i$.
$F_{312}(x,q,z)$ is given by the following continued fraction. $$
F_{312}(x,q,z)=\frac{1}{1-z K_0+\cfrac{z}{1-z K_1+\cfrac{z}{1-z
K_2+ \cfrac{z}{1-z K_3+\cfrac{z}{\ddots}}}}},$$ where
$K_n=(x-1)\C_n q^n z^n + (q-1)\C_{<n}(qz)$ for $n\ge0$.
\end{theorem}

Note that the first values of $K_n$ are \bce\bt{ll} $K_0=x-1$, &
$K_1=(x-1)qz+q-1$, \\ $K_2=2(x-1)q^2 z^2+(q-1)(1+qz)$, &
$K_3=5(x-1)q^3 z^3+(q-1)(1+qz+2q^2z^2)$.\et\ece

\begin{proof}
We use the same decomposition as above, now keeping track of
tunnels of depth $<0$ as well. For every $h\ge0$, let
$\tdneg^h(D)$ be the number of tunnels of $D$ of height $h$ and
length $<2(h+1)$. Note that $\tdneg(D)=\sum_{h\ge0}\tdneg^h(D)$.
To follow the same structure as in the previous proof, counting
tunnels height by height, it will be convenient that at the $h$-th
step of the iteration, $q$ marks not only tunnels of depth $<0$ up
to height $h$ but also all the tunnels at higher levels. Denote by
$\at^{>h}(D)$ the number of tunnels of $D$ of height strictly
greater than $h$.

We will show now that for every $h\ge1$, the generating function
for Dyck paths where $x$ marks the statistic
$\tdzero^0(\cdot)+\cdots+\tdzero^{h-1}(\cdot)$ and $q$ marks
$\tdneg^0(\cdot)+\cdots+\tdneg^{h-1}(\cdot)+\at^{>h-1}(\cdot)$ is given by the
continued fraction of the theorem truncated at level $h$, with the
$(h+1)$-st level replaced with $\C(qz)$.

The analogous to equation (\ref{height0}) is now
\beq\sum_{D\in\D}x^{\tdzero^0(D)}q^{\tdneg^0(D)+\at^{>0}(D)}z^{|D|}=\frac{1}{1-z[x-1+\C(qz)]}.\label{hei0}\eeq
Indeed, decomposing $D$ as $\u A_1\d\cdots\u A_r\d$, $q$ counts
all the tunnels that appear in any $A_i$, and whenever an $A_i$ is
empty we must mark it as $x$.

Let us enumerate now tunnels of depth 0 and depth $<0$ at both
height~0 and height~1. Modifying (\ref{hei0}) so that $q$ no
longer counts tunnels at height~1, we get
\beq\sum_{D\in\D}x^{\tdzero^0(D)}q^{\tdneg^0(D)+\at^{>1}(D)}z^{|D|}=\frac{1}{1-z\left[x-1+\dfrac{1}{1-z
\C(qz)}\right]},\label{hei0mod}\eeq which corresponds to writing
each $A_i$ as $A_i=\u B_1\d\cdots\u B_s\d$, and having $q$ count
all tunnels in each $B_j$. Now, in order for $x$ to mark tunnels
of depth 0 at height~1, each $B_j=\u\d$, that in (\ref{hei0mod})
is counted as $qz$, has to be now counted as $xqz$ instead.
Similarly, to have $q$ mark tunnels of depth $<0$ at height~1, we
must count each empty $B_j$ as $q$, not as 1. This gives us the
following GF:
\bea\nn\sum_{D\in\D}x^{\tdzero^0(D)+\tdzero^1(D)}q^{\tdneg^0(D)+\tdneg^1(D)+\at^{>1}(D)}z^{|D|}\hspace{2cm}\\
\nn =
\frac{1}{1-z\left[x-1+\dfrac{1}{1-z[(x-1)qz+q-1+\C(qz)}\right]}.\eea

Iterating this process level by level indefinitely we obtain the
continued fraction of the theorem. At each step, from the GF where
$x$ marks $\tdzero^0(\cdot)+\cdots+\tdzero^{h-1}(\cdot)$, and $q$ marks
$\tdneg^0(\cdot)+\cdots+\tdneg^{h-1}(\cdot)+\at^{>h-1}(\cdot)$, we can obtain
the one where $x$ marks $\tdzero^0(\cdot)+\cdots+\tdzero^h(\cdot)$ and $q$
marks $\tdneg^0(\cdot)+\cdots+\tdneg^h(\cdot)+\at^{>h}(\cdot)$ by replacing
the $\C(qz)$ at the lowest level with \beq\label{hthstep}
\frac{1}{1-z[(x-1)\C_h q^h z^h+(q-1)\C_{<h}(qz)+\C(qz)]}.\eeq This
change makes $x$ account for tunnels of depth 0 at height $h$,
which in the decomposition correspond to the $\C_h$ possible
elevated Dyck paths of length $2(h+1)$ when they occur at height
$h$. It also makes $q$ count tunnels of depth $<0$ at height $h$,
which in the decomposition correspond to elevated Dyck paths at
height $h$ of length $<2(h+1)$. The GF for these ones becomes $q\C_{<h}(qz)$
instead of $\C_{<h}(qz)$, since for every $j<h$, an elevated path
$\u C\d$ with $C\in\D_j$ contributes to one extra tunnel of depth
$<0$ at height $h$, aside from the $j$ tunnels of height $>h$ that
it contains.
\end{proof}

For 231-avoiding permutations we get the following GF.

\begin{corollary} Let $\C_{<i}(z)$ be as in Theorem~\ref{th:312}.
$F_{231}(x,q,z)$ is given by the following continued fraction.
$$ F_{231}(x,q,z)=\frac{1}{1-zK'_0+\cfrac{qz}{1-zK'_1+\cfrac{qz}{1-zK'_2+\cfrac{qz}{1-zK'_3
+\cfrac{qz}{\ddots}}}}},$$ where $K'_n=(x-q)\C_n z^n +
(1-q)\C_{<n}(z)$.
\end{corollary}

The first values of $K_n$ are \bce\bt{ll} $K'_0=x-q$, &
$K'_1=(x-q)z+1-q$, \\ $K'_2=2(x-q)z^2+(1-q)(1+z)$, &
$K'_3=5(x-q)z^3+(1-q)(1+z+2z^2)$.\et\ece

\begin{proof}
By Lemma~\ref{lemma:sim}, we have that
$F_{231}(x,q,z)=F_{312}(\frac{x}{q},\frac{1}{q},qz),$ so the
expression follows from Theorem~\ref{th:312}.
\end{proof}

\section{Double restrictions}\label{sec:double}

In this section we consider simultaneous avoidance of any two
patterns of length~3. Using Lemma~\ref{lemma:sim}, the pairs of
patterns fall into the following equivalence classes.

\bce{\bf a) $\{123,132\}\simfe\{123,213\}$ \\
b) $\{231,321\}\quad\simf\quad$ b') $\{312,321\}$ \\
c) $\{132,213\}$ \\ d) $\{231,312\}$ \\
e) $\{132,231\}\simfe\{213,231\}\quad\simf\quad$
e') $\{132,312\}\simfe\{213,312\}$ \\
f) $\{132,321\}\simfe\{213,321\}$ \\
g) $\{123,231\}\quad\simf\quad$
g') $\{123,312\}$ \\
h) $\{123,321\}$
 }\ece

In \cite{SS85} it is shown that the number of permutations in
$S_n$ avoiding any of the pairs in the classes {\bf a)}, {\bf b)},
{\bf b')}, {\bf c)}, {\bf d)}, {\bf e)}, and {\bf e')} is
$2^{n-1}$, and that for the pairs in {\bf f)}, {\bf g)} and {\bf
g')}, the number of permutations avoiding any of them is ${n
\choose 2}+1$. The case {\bf h)} is trivial because this pair is
avoided only by permutations of length at most 4.

In terms of generating functions, this means that when we
substitute $x=q=1$ in $F_\Sigma(x,q,z)$, where $\Sigma$ is any of
the pairs in the classes from {\bf a)} to {\bf e')}, we get
$$F_\Sigma(1,1,z)=\sum_{n\ge0}2^{n-1}z^n=\frac{1-z}{1-2z}.$$
If $\Sigma$ is a pair from the classes {\bf f)}, {\bf g)}, {\bf
g')}, we get
$$F_\Sigma(1,1,z)=\sum_{n\ge0}({n \choose 2}+1)z^n=\frac{1-2z+2z^2}{(1-z)^3}.$$

\subsection{a) $\{123,132\}\simfe\{123,213\}$}

\begin{prop}\label{th:123132}
\bea\nn F_{\{123,132\}}(x,q,z)=F_{\{123,213\}}(x,q,z)\hspace{8cm} \\
\nn
\hspace*{1cm}=\frac{1+xz+(x^2-4q)z^2+(-3xq+q+q^2)z^3+(xq+xq^2-3x^2q+3q^2)z^4}{(1-qz^2)(1-4qz^2)}.\eea
\end{prop}

\begin{proof}
Consider the bijection $\kra:\S_n(132)\longrightarrow\D_n$
described in Subsection~\ref{sec:bij}. It is shown in \cite{Kra01}
that the height of the Dyck path $\kra(\pi)$ is the length of the
longest increasing subsequence of $\pi$. In particular,
$\pi\in\S_n(12\cdots(k+1),132)$ if and only if $\kra(\pi)$ has
height at most $k$. Thus, by Lemma~\ref{lemma:kra},
$F_{\{123,132\}}(x,q,z)$ can be written in terms of Dyck paths as
\beq\label{F2aDyck} \sum_{D\in\D^{\le
2}}x^{\ct(D)}q^{\rt(D)}z^{|D|}.\eeq

Let us first find the univariate GF for paths of height at most 2 (with no
statistics). Clearly, the GF for Dyck paths of height at most 1 is
$\frac{1}{1-z}$, since such paths are just sequences of hills. A
path $D$ of height at most 2 can be written uniquely as $D=\u
A_1\d\u A_2\d\cdots\u A_r\d$, where each $A_i$ is a path of height
at most 1. The GF for each $\u A_i\d$ is $\frac{z}{1-z}$. Thus,
$$\sum_{D\in\D^{\le 2}}z^{|D|}=\frac{1}{1-\frac{z}{1-z}}=\frac{1-z}{1-2z}=\sum_{n\ge0}2^{n-1}z^n.$$

In the rest of this proof, we assume that all Dyck paths that
appear have height at most 2 unless otherwise stated. To compute
(\ref{F2aDyck}), we will separate paths according to their height
at the middle. Consider first paths whose height at the middle is
0. Splitting such a path at its midpoint we obtain a pair of paths
of the same length. Thus, the corresponding GF is
\beq\label{F2a_h0}\sum_{m\ge0}2^{m-1}z^m\cdot2^{m-1}q^m
z^m=\frac{1-3qz^2}{1-4qz^2},\eeq since the number of right tunnels
of such a path is the semilength of its right half.

\begin{figure}[hbt]
\epsfig{file=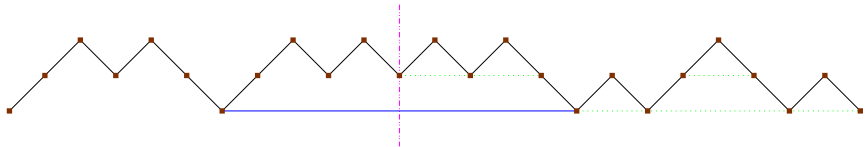,width=4.5in} \caption{\label{fig:height2}
A path of height 2 with a centered tunnel.}
\end{figure}

Now we consider paths whose height at the middle is 1. It is easy
to check that without the variables $x$ and $q$, the GF for such
paths is \beq\label{G21}\frac{z}{1-4z^2}.\eeq Let us first look at
paths of this kind that have a centered tunnel. They must be of
the form $D=A\u B\d C$ where $A,C\in\D^{\le 2}$ have the same
length, and $B$ is a sequence of an even number of hills. Thus,
their GF is
\beq\label{eqct}xz\cdot\frac{1}{1-qz^2}\cdot\frac{1-3qz^2}{1-4qz^2},\eeq
where $x$ marks the centered tunnel, $\frac{1}{1-qz^2}$
corresponds to the sequence of hills $B$, half of which give right
tunnels, and the last fraction comes from the pair $AC$, which is
counted as in (\ref{F2a_h0}). From (\ref{G21}) and (\ref{eqct}) it
follows that the univariate GF (with just variable $z$) for paths
with height at the middle 1, not having a centered tunnel, is
$$\frac{z}{1-4z^2}-\frac{z(1-3z^2)}{(1-z^2)(1-4z^2)}=\frac{2z^3}{(1-z^2)(1-4z^2)}.$$
By symmetry, in half of these paths, the tunnel of height~0 that
goes across the middle is a right tunnel. Thus, the multivariate
GF for all paths with height~1 at the middle is
\beq\label{F2a_h1}\frac{xz(1-3qz^2)}{(1-qz^2)(1-4qz^2)}+\frac{(q+1)qz^3}{(1-qz^2)(1-4qz^2)}.\eeq
Here the right summand corresponds to paths with no centered
tunnel: the term $(q+1)$ distinguishes whether the tunnel that
goes across the middle is a right tunnel or not, and the other
$q$'s mark tunnels completely contained in the right half.

Paths with height 2 at the middle are easy to enumerate now.
Indeed, they must have a peak $\u\d$ in the middle, whose removal
induces a bijection between these paths and paths with height~1 at
the middle. This bijection preserves the number of right tunnels,
and decreases the length and the number of centered tunnels by
one. Thus, the GF for paths with height 2 at the middle is $xz$
times expression (\ref{F2a_h1}). Adding up this GF, for paths with
height 2 at the middle, to the expressions (\ref{F2a_h0}) and
(\ref{F2a_h1}), for paths with height at the middle 0 and 1
respectively, we obtain the desired expression for
$F_{\{123,132\}}(x,q,z)$.
\end{proof}

Let us see how the same technique used in this proof can be generalized to
enumerate fixed points in $\S_n(132,12\cdots(k+1))$ for an
arbitrary $k\ge0$.

\begin{theorem}
For $k\ge 0$, let
$$M_k(x,z):=F_{\{132,12\cdots(k+1)\}}(x,1,z)=\sum_{n\ge
0}\sum_{\pi\in\S_n(132,12\cdots(k+1))} x^{\fp(\pi)} z^n.$$ Then
the $M_k$'s satisfy the recurrence: \beq M_k(x,z)=\sum_{\ell=0}^k
G_{k,\ell}(z)(1+(x-1)z M_{\ell-1}(x,z)),\eeq where
$M_{-1}(x,z):=0$, and $G_{k,\ell}(z)$ is defined as
$$G_{k,\ell}(z):=\sum_{n\ge 0}g_{k,\ell}^2(n) z^n,\mbox{ where }
\sum_{n\ge 0}g_{k,\ell}(n) z^n=\frac{U_{\ell}(\frac{1}{2z})}{z
U_{k+1}(\frac{1}{2z})},$$ where $U_m$ are the Chebyshev
polynomials of the second kind, defined by the recurrence
$$U_0(t)=1,\ U_1(t)=2t,\ \mbox{and}\
U_r(t)=2tU_{r-1}(t)-U_{r-2}(t).$$
\end{theorem}

Before proving this theorem, let us show how to apply it to
obtain the GFs $M_k$ for the first few values of $k$. For $k=1$, we have
$G_{1,0}(z)=\frac{z}{1-z^2}$, $G_{1,1}(z)=\frac{1}{1-z^2}$, so
$$M_1(x,z)=\frac{1+xz}{1-z^2}.$$

For $k=2$, we get $G_{2,0}(z)=\frac{z^2}{1-4z^2}$,
$G_{2,1}(z)=\frac{z}{1-4z^2}$, $G_{2,2}(z)=1+\frac{z^2}{1-4z^2}$,
so
$$M_2(x,z)=\frac{1+xz+(x^2-4)z^2+(2-3x)z^3+(3+2x-3x^2)z^4}{(1-z^2)(1-4z^2)},$$
which is the expression of Proposition~\ref{th:123132} for $q=1$.

For $k=3$, we obtain
$G_{3,0}(z)=\frac{z^3+z^5}{(1-z^2)(1-7z^2+z^4)}$,
$G_{3,1}(z)=\frac{z^2+z^4}{(1-z^2)(1-7z^2+z^4)}$,
$G_{3,2}(z)=\frac{z(1-4z^2+z^4)}{(1-z^2)(1-7z^2+z^4)}$,
$G_{3,2}(z)=1+\frac{z^2(1-4z^2+z^4)}{(1-z^2)(1-7z^2+z^4)}$, so \\
$M_3(x,z)=$ {\small
$[1+xz+(x^2-12)z^2+(x^3-11x+2)z^3+(-10x^2+4x+45)z^4+(-10x^3+4x^2+37x-10)z^5$
$\hspace*{17mm}+(25x^2-22x-52)z^6
+(25x^3-22x^2-41x+16)z^7+(-12x^2+16x+16)z^8+(-12x^3+16x^2$
$\hspace*{17mm}+12x-8)z^9] \ / \
[(1-z^2)^2(1-4z^2)(1-7z^2+z^4)]$}.

\begin{proof}
As shown at the beginning of the previous proof, $\kra$ induces a
bijection between $\S_n(132,12\cdots(k+1))$ and $\D^{\le k}$, the
set of Dyck paths of height at most $k$. Thus, by
Lemma~\ref{lemma:kra}, $$M_k(x,z)=\sum_{D\in\D^{\le
k}}x^{\ct(D)}z^{|D|}.$$

In order to find a recursion for this GF, we are going to count
pairs $(D,S)$ where $D\in\D^{\le k}$ and $S$ is a subset of
$\mathrm{CT}(D)$, the set of centered tunnels of $D$. In other
words, we are considering Dyck paths where some centered tunnels
(namely those in $S$) are marked. Each such pair is given weight
$(x-1)^{|S|} z^{|D|}$, so that for a fixed $D$, the sum of weights
of all pairs $(D,S)$ is $\sum_{S\subseteq
\mathrm{CT}(D)}{(x-1)^{|S|} z^{|D|}}=((x-1)+1)^{|\mathrm{CT}(D)|}
z^{|D|}=x^{\ct(D)} z^{|D|}$, which is precisely the weight that
$D$ has in $M_k(x,z)$.

\ms

\begin{figure}[hbt]
\epsfig{file=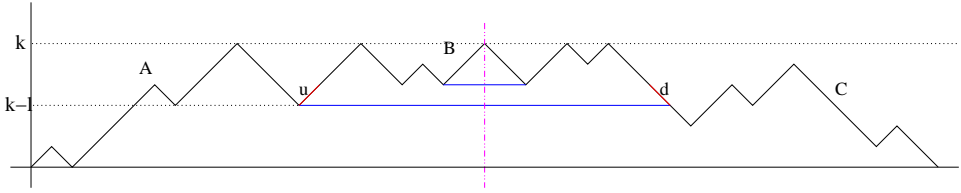,width=5in} \caption{\label{fig:heightk} A
path of height $k$ with two marked centered tunnels.}
\end{figure}

If $D\in\D^{\le k}$ has some marked centered tunnel, consider the
decomposition $D=A\u B\d C$ given by the longest marked tunnel
(i.e., all the other marked tunnels are inside the part $B$ of the
path). Let $\ell$ be the distance between this tunnel and the line
$y=k$ (see Figure~\ref{fig:heightk}). Equivalently, $A$ ends at height $k-\ell$, the same height
where $C$ begins. Then, $B$ is an arbitrary Dyck path of height at
most $\ell-1$ with possibly some marked centered tunnels, so its
corresponding GF is $M_{\ell-1}(x,z)$ (with the convention
$M_{-1}(x,z):=0$, since for $\ell=0$ there is no such $B$). Giving
weight $(x-1)$ to the tunnel that determines our decomposition, we
have that the part $\u B\d$ of the path contributes
$(x-1)zM_{\ell-1}(x,z)$ to the GF.

Now we look at the GF for the part $A$ of the path. Let
$g_{k,\ell}(n)$ be the number of paths from $(0,0)$ to
$(n,k-\ell)$ staying always between $y=0$ and $y=k$. A path of
this type can be decomposed uniquely as $A=E_k\u E_{k-1}\u \cdots
\u E_{\ell+1} \u E_\ell$, where each $E_i\in\D^{\le i}$. The GF of
Dyck paths of height at most $i$ is
$$M_i(1,z)=\frac{U_i(\frac{1}{2\sqrt{z}})}{\sqrt{z}U_{i+1}(\frac{1}{2\sqrt{z}})},$$
as shown for example in \cite{Kra01}. Let $w=\sqrt{z}$, which is
the weight of a single step of a path, and let
$\wt{G}_{k,\ell}(w):=\sum_{n\ge 0}g_{k,\ell}(n) w^n$. From the
above decomposition of $A$,
$$\wt{G}_{k,\ell}(w)=M_k(1,w^2)w M_{k-1}(1,w^2)w \cdots w
M_\ell(1,w^2)=\frac{U_\ell(\frac{1}{2 w})}{w
U_{k+1}(\frac{1}{2w})}.$$ The part $C$ of the path $D$, flipped
over a vertical line, can be regarded as a path with the same
endpoints as $A$, since it must have the same length and end at
the same height $k-\ell$. Thus, the GF for pairs $(A,C)$ of paths
of the same length from height~0 to height $k-\ell$ and not going
above $y=k$ is $\sum_{n\ge 0}g_{k,\ell}^2(n) z^n=G_{k,\ell}(z)$.

Hence, the GF for paths $D\in\D^{\le k}$ having the longest marked
centered tunnel at height $k-\ell$ is
$G_{k,\ell}(z)(x-1)zM_{\ell-1}(x,z)$.

\ms

If $D$ has no marked tunnel, decompose it as $D=AC$ where $A$ and
$C$ have the same length. Letting $k-\ell$ be again the height
where $A$ ends and $C$ begins, the situation is the same as above
but without any contribution coming from the central part of $D$.
The parameter $\ell$ can take any value between $0$ and $k$. Thus,
summing over all possible decompositions of $D$, we get
$$M_k(x,z)=\sum_{\ell=0}^k
G_{k,\ell}(z)(1+(x-1)z M_{\ell-1}(x,z)).$$
\end{proof}

\subsection{b, b') $\{231,321\}\simf\{312,321\}$}

\begin{prop}\label{th:312321}
\beq\label{eq:312321} F_{\{312,321\}}(x,q,z)
=\frac{1-qz}{1-(x+q)z+(x-1)qz^2}.\eeq
\end{prop}

\begin{proof}
The length of the longest decreasing subsequence of $\pi$ equals
the height of the Dyck path $\kra(\bar\pi)$. In particular, we
have a bijection \bce\bt{ccc} $\S_n(312,321)$ &
$\longleftrightarrow$ & $\D^{\le 2}_n$ \\ $\pi$ & $\mapsto$ &
$\kra(\bar\pi)$ \et\ece

Thus, by Lemma~\ref{lemma:kra},
$$F_{\{312,321\}}(x,q,z)=\sum_{D\in\D^{\le 2}}x^{\tdzero(D)}q^{\tdneg(D)}z^{|D|}.$$
But the only tunnels of depth 0 that a Dyck path of height at most
2 can have are hills, and the only tunnels of depth $<0$ that it
can have are peaks at height 2. A path $D\in\D^{\le 2}$ can be
written uniquely as $D=\u A_1\d\u A_2\d\cdots\u A_r\d$, where each
$A_i$ is a (possibly empty) sequence of hills. An empty $A_i$
creates a tunnel of depth 0 in $D$, so it contributes as $x$. An
$A_i$ of length $2j>0$ contributes as $q^jz^j$, since it creates
$j$ peaks at height 2 in $D$. Thus,
$$F_{\{312,321\}}(x,q,z)=\frac{1}{1-z\left(x+\dfrac{qz}{1-qz}\right)},$$
which is equivalent to (\ref{eq:312321}).
\end{proof}

\begin{corollary}
$$F_{\{231,321\}}(x,q,z) =\frac{1-z}{1-(x+1)z+(x-q)z^2}.$$
\end{corollary}

\begin{proof}
It follows from Lemma~\ref{lemma:sim} that
$F_{\{231,321\}}(x,q,z)=F_{\{312,321\}}(\frac{x}{q},\frac{1}{q},qz)$.
\end{proof}

As in the previous section, these results can be generalized to
the case when instead of the pattern 321 we have a decreasing
pattern $(k+1)k\cdots21$ of arbitrary length. For $i,h\ge 0$, let
$\C_i^{\le h}$ be the number of Dyck paths of length $2i$ and
height at most $h$. As mentioned before, it is known that
$$\sum_{i\ge 0}\C_i^{\le
h}z^i=\frac{U_h(\frac{1}{2\sqrt{z}})}{\sqrt{z}U_{h+1}(\frac{1}{2\sqrt{z}})},$$
where $U_m$ are the Chebyshev polynomials of the second kind. Let
$$\C_{<i}^{\le h}(z):=\sum_{j=0}^{i-1}\C_j^{\le h} z^j.$$
The following theorem deals with fixed points and excedances in
$\S_n(312,(k+1)k\cdots1)$ for any $k\ge0$.

\begin{theorem}\label{th:312k1}
Let $\C_i^{\le h}=|\D_i^{\le h}|$ and $\C_{<i}^{\le h}(z)$ be
defined as above. Then, for $k\ge0$,
$$F_{\{312,(k+1)k\cdots1\}}(x,q,z)=A_0^k(x,q,z),$$ where $A_i^k$ is recursively defined
by $$ A_i^k(x,q,z)=\left\{ \bt{ll} $\dfrac{1}{1-z[(x-1)\C_i^{\le
k-i-1}q^iz^i+(q-1)\C_{<i}^{\le
k-i-1}(qz)+A_{i+1}^k(x,q,z)]}$ & if $i<k$, \\
$1$ & if $i=k$. \et\right.$$
\end{theorem}

For example, for $k=2$ we obtain Proposition~\ref{th:312321}, and
for $k=3$ we get \bea\nn
F_{\{312,4321\}}(x,q,z)&=&\frac{1}{1-z\left[x-1+\dfrac{1}{1-z\left[(x-1)qz+q-1+\frac{1}{1-qz}\right]}\right]} \\
\nn &=&
\frac{1-2qz+(q^2-xq)z^2+(xq^2-q^2)z^3}{1-(x+2q)z+(xq+q^2-q)z^2+(x^2q-xq)z^3+(-x^2q^2+2xq^2-q^2)z^4}.\eea

\begin{proof}
It is analogous to the proof of Theorem~\ref{th:312}, with the
only difference that here we consider only those paths that do not
go above the line $y=k$.
\end{proof}

Making the appropriate substitutions in the statement of
Theorem~\ref{th:312k1}, we obtain an expression for the generating
function
$F_{\{231,(k+1)k\cdots1\}}(x,q,z)=F_{\{312,(k+1)k\cdots1\}}(\frac{x}{q},\frac{1}{q},qz)$.

\subsection{c) $\{132,213\}$}

\begin{prop}\label{th:132213}
$$F_{\{132,213\}}(x,q,z)
=\frac{1-(1+q)z-2qz^2+4q(1+q)z^3-(xq^2+xq+5q^2)z^4+2xq^2z^5}{(1-z)(1-xz)(1-qz)(1-4qz^2)}.$$
\end{prop}

\begin{proof}
We use again the bijection $\kra:\S_n(132)\longrightarrow\D_n$.
From its description given in Subsection~\ref{sec:bij}, it is not
hard to see that a permutation $\pi\in\S_n(132)$ avoids 213 if and
only if all the valleys of the corresponding Dyck path $\kra(\pi)$
have their lowest point on the $x$-axis. A path with such property
can be described equivalently as a \emph{sequence of pyramids}.
Denote by $\Pyr_n\subseteq\D_n$ the set of sequences of pyramids
of length $2n$, and let $\Pyr:=\bigcup_{n\ge0}\Pyr_n$. We have
just seen that $\kra$ restricts to a bijection between
$\S_n(132,213)$ and $\Pyr_n$. By Lemma~\ref{lemma:kra}, we can
write $F_{\{132,213\}}(x,q,z)$ as
$$\sum_{D\in\Pyr} x^{\ct(D)}q^{\rt(D)}z^{|D|}.$$
Since for each $n\ge1$ there is exactly one pyramid of length
$2n$, the univariate GF of sequences of pyramids is just
$\sum_{D\in\Pyr}z^{|D|}=\frac{1}{1-\frac{z}{1-z}}=\frac{1-z}{1-2z}=1+\sum_{n\ge1}2^{n-1}z^n$.

Let us first consider elements of $\Pyr$ that have height~0 in the
middle (equivalently, the two central steps are $\d\u$). Each one
of their halves is a sequence of pyramids, both of the same
length. They have no centered tunnels, and the number of right
tunnels is given by the semilength of the right half. Thus, their
multivariate GF is \beq\label{pyr1} 1+\sum_{m\ge1}4^{m-1}q^m
z^{2m}=1+\frac{qz^2}{1-4qz^2}.\eeq

Now we count elements of $\Pyr$ whose two central steps are
$\u\d$. They are obtained uniquely by inserting a pyramid of
arbitrary length in the middle of a path with height~0 at the
middle. The tunnels created by the inserted pyramid are all
centered tunnels, so the corresponding GF is
\beq\label{pyr2}\frac{xz}{1-xz}\left(1+\frac{qz^2}{1-4qz^2}\right).\eeq

\begin{figure}[hbt]
\epsfig{file=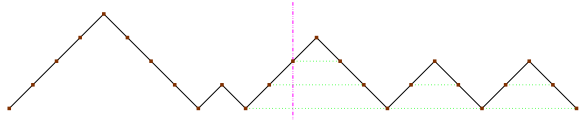,width=3in} \caption{\label{fig:seqpyr} A
sequence of pyramids.}
\end{figure}

It remains to count the elements of $\Pyr$ that in the middle have
neither a peak nor a valley. From a non-empty sequence of pyramids
with height~0 in the middle, if we increase the size of the
leftmost pyramid of the right half by an arbitrary number of
steps, we obtain a sequence of pyramids whose two central steps
are $\u\u$. Reciprocally, by this procedure every such sequence of
pyramids can be obtained in a unique way from a sequence of
pyramids with height~0 in the middle. Thus, the GF for the
elements of $\Pyr$ whose two central steps are $\u\u$ is
\beq\label{pyr3}\frac{qz}{1-qz}\cdot\frac{qz^2}{1-4qz^2}.\eeq By
symmetry, the GF for the elements of $\Pyr$ whose two central
steps are $\d\d$ is
\beq\label{pyr4}\frac{z}{1-z}\cdot\frac{qz^2}{1-4qz^2},\eeq where
the difference with respect to~\ref{pyr3} is that now the pyramid
across the middle does not create right tunnels. Adding up
(\ref{pyr1}), (\ref{pyr2}), (\ref{pyr3}) and (\ref{pyr4}) we get
the desired GF.
\end{proof}

\subsection{d) $\{231,312\}$}

\begin{prop}\label{th:231312}
$$F_{\{231,312\}}(x,q,z)
=\frac{1-qz^2}{1-xz-2qz^2}.$$
\end{prop}

\begin{proof}
We have shown in the proof of Proposition~\ref{th:132213} that
$\kra$ induces a bijection between $\S_n(132,213)$ and $\Pyr_n$,
the set of sequences of pyramids of length $2n$. Composing it with
the complementation operation, we get a bijection
$\pi\mapsto\kra(\bar\pi)$ between $\S_n(231,312)$ and $\Pyr_n$.
Together with Lemma~\ref{lemma:kra}, this allows us to express
$F_{\{231,312\}}(x,q,z)$ as
$$\sum_{D\in\Pyr} x^{\tdzero(D)}q^{\tdneg(D)}z^{|D|}.$$
All that remains is to observe how many tunnels of zero and
negative depth are created by a pyramid according to its size. A
pyramid of odd semilength $2m+1$ creates one tunnel of depth 0 and
$m$ tunnels of depth $<0$. A pyramid of even semilength $2m$
creates only $m$ tunnels of depth $<0$. Thus, we have that
$$F_{\{231,312\}}(x,q,z)=\frac{1}{1-\dfrac{xz}{1-qz^2}-\dfrac{qz^2}{1-qz^2}},$$
which equals the expression above.
\end{proof}

\subsection{e, e') $\{132,231\}\simfe\{213,231\}\simf\{132,312\}\simfe\{213,312\}$}

\begin{prop}\label{th:132231}
$$F_{\{132,231\}}(x,q,z)=F_{\{213,231\}}(x,q,z)
=\frac{1-z-qz^2+xqz^3}{(1-xz)(1-z-2qz^2)}.$$
\end{prop}

\begin{proof}
As usual, we use the bijection
$\kra:\S_n(132)\longrightarrow\D_n$. Now we are interested in how
the condition that $\pi$ avoids 231 is reflected in the Dyck path
$\kra(\pi)$. It is easy to see from the description of $\kra$ and
$\kra^{-1}$ in Subsection~\ref{sec:bij} that $\pi$ is 231-avoiding
if and only if $\kra(\pi)$ does not have any two consecutive
up-steps after the first down-step (equivalently, all the
non-isolated up-steps occur at the beginning of the path). Let
$\E_n\subseteq\D_n$ be the set of Dyck paths with this condition,
and let $\E:=\bigcup_{n\ge0}\E_n$. Then, $\kra$ induces a
bijection between $S_n(132,231)$ and $\E_n$. By
Lemma~\ref{lemma:kra}, $F_{\{132,231\}}(x,q,z)$ can be written as
$$\sum_{D\in\E} x^{\ct(D)}q^{\rt(D)}z^{|D|}.$$

\begin{figure}[hbt]
\epsfig{file=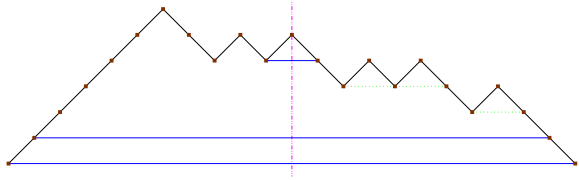,width=3in} \caption{\label{fig:roller} A
path in $\E$ with a peak in the middle and two bottom tunnels.}
\end{figure}

If $D\in\E$, centered tunnels of $D$ can appear only in the
following two places. There can be a centered tunnel produced by a
peak in the middle of $D$. All the other centered tunnels of $D$
must have their endpoints in the initial ascending run and the
final descending run of $D$ (that is, in their corresponding
decomposition $D=A\u B\d C$, $A$ is a sequence of up-steps and $C$
is a sequence of down-steps). For convenience we call this second
kind of tunnels \emph{bottom} tunnels. All the right tunnels of
$D$ come from peaks on the right half.

It is an exercise to check that the number of paths in $\E_n$
having a peak in the middle and $r$ peaks on the right half is
$\binom{n-r-1}{r}2^{r-1}$ if $r\ge1$, and 1 if $r=0$. Similarly,
the number of paths in $\E_n$ with no peak in the middle and $r$
peaks on the right half is $\binom{n-r}{r}2^{r-1}$ if $r\ge1$, and
0 if $r=0$. Let us ignore for the moment the bottom tunnels.
For peaks in the middle and right tunnels we have the following
GF. \bea\label{eqP} P(x,q,z):=\sum_{D\in\E\setminus\E_0}
x^{\#\{\mathrm{peaks\ in\ the\ middle\ of\ D}\}} q^{\rt(D)}z^{|D|}\hspace{4cm}\\
\nn=\sum_{n\ge1}\left[\sum_{r=1}^{\lfloor
n/2\rfloor}\binom{n-r}{r}2^{r-1}q^r+x\left( 1+\sum_{r=1}^{\lfloor
(n-1)/2\rfloor}\binom{n-r-1}{r}2^{r-1}q^r\right)\right]z^n
=\frac{xz+(q-x)z^2-xqz^3}{(1-z)(1-z-2qz^2)}.\hspace{-6mm}\eea Now,
to take into account all centered tunnels, we use that every
$D\in\E$ can be written uniquely as $D=\u^k D'\d^k$, where $k\ge0$
and $D'\in\E$ has no bottom tunnels. The GF for elements of $\E$
that do have bottom tunnels, where $x$ marks peaks in the middle,
is $zx+zP(x,q,z)$ (the term $zx$ is the contribution of the path
$\u\d$). Hence, the sought GF where $x$ marks all centered tunnels
is
$$F_{\{132,231\}}(x,q,z)=\frac{1}{1-xz}[1+P(x,q,z)-zx-zP(x,q,z)]=1+\frac{1-z}{1-xz}P(x,q,z),$$
which together with (\ref{eqP}) implies the proposition.
\end{proof}

\begin{corollary}
$$F_{\{132,312\}}(x,q,z)=F_{\{213,312\}}(x,q,z)=\frac{1-qz-qz^2+xqz^3}{(1-xz)(1-qz-2qz^2)}.$$
\end{corollary}

\begin{proof}
It follows from Lemma~\ref{lemma:sim} that
$F_{\{132,312\}}(x,q,z)=F_{\{132,231\}}(\frac{x}{q},\frac{1}{q},qz)$.
\end{proof}

\subsection{f) $\{132,321\}\simfe\{213,321\}$}

\begin{prop}\label{th:132321}
$$F_{\{132,321\}}(x,q,z)=F_{\{213,321\}}(x,q,z)
=\frac{1-(1+q)z+2qz^2}{(1-z)(1-xz)(1-qz)}.$$
\end{prop}

\begin{figure}[hbt]
\epsfig{file=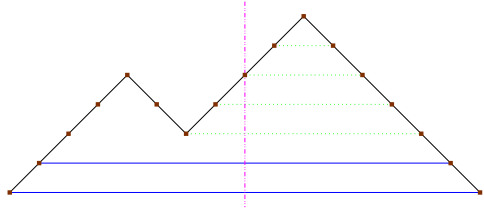,width=2.5in}
\caption{\label{fig:twopeaks} A path with two peaks.}
\end{figure}

\begin{proof}
From the definition of the bijection
$\kra:\S_n(132)\longrightarrow\D_n$ and the description of its
inverse given in Subsection~\ref{sec:bij}, it follows that the
number of peaks of the Dyck path $\kra(\pi)$ equals the length of
the longest decreasing subsequence of $\pi$. In particular, $\pi$
is 321-avoiding if and only if $\kra(\pi)$ has at most two peaks.
By Lemma~\ref{lemma:kra}, $F_{\{132,321\}}(x,q,z)=\sum
x^{\ct(D)}q^{\rt(D)}z^{|D|}$, where the sum is over Dyck paths $D$
with at most two peaks. Clearly, such a path can be uniquely
written as $D=\u^k D'\d^k$, where $k\ge0$ and $D'$ is either empty
or a pair of adjacent pyramids (see Figure~\ref{fig:twopeaks}). Therefore,
$$F_{\{132,321\}}(x,q,z)=\frac{1}{1-xz}\left(1+\frac{z}{1-z}\cdot\frac{qz}{1-qz}\right),$$
since centered tunnels are produced by the steps outside $D'$, and
right tunnels are created by the right pyramid of $D'$.
\end{proof}

This case can be generalized to the situation when instead of 321
we have a decreasing pattern of arbitrary length. Observe that by
Lemma~\ref{lemma:sim},
$F_{\{132,(k+1)k\cdots21\}}(x,q,z)=F_{\{213,(k+1)k\cdots21\}}(x,q,z)$
for all $k$.

\begin{theorem}
$$\sum_{k\ge0}F_{\{132,(k+1)k\cdots21\}}(x,q,z)\, p^k
=\frac{2(1+xz(p-1))}{(1-p)[1+(1+q-2x)z-qz^2
(p-1)^2+\sqrt{f_1(q,z)}]},$$ where $f_1(q,z)=
1-2(1+q)z+[(1-q)^2-2q(p-1)(p+3)]z^2-2q(1+q)(p-1)^2 z^3+q^2 (p-1)^4
z^4$.
\end{theorem}

\begin{proof}
As in the previous proof, we use the fact that the number of peaks
of $\kra(\pi)$ equals the length of the longest decreasing
subsequence of $\pi$. Thus, $\kra$ induces a bijection between
$S_n(132,(k+1)k\cdots21)$ and the subset of $\D_n$ of paths with
at most $k$ peaks. This implies that
$\sum_{k\ge0}F_{\{132,(k+1)k\cdots21\}}(x,q,z)\, p^k$ can be
expressed as $$\frac{1}{1-p}\sum_{D\in\D} x^{\ct(D)} q^{\rt(D)}
p^{\#\{\mathrm{peaks\ of\ }D \}} z^{|D|}.$$ The result now follows
from \cite[Theorem 6.5]{EliDeu} and the expression for $\sum
x^{\ct(D)} q^{\rt(D)} p^{\#\{\mathrm{peaks\ of\ }D \}} z^{|D|}$
given in its proof.
\end{proof}

\subsection{g, g') $\{123,231\}\simf\{123,312\}$}

\begin{prop}\label{th:123231}
$$F_{\{123,312\}}(x,q,z)
=\frac{1+xz+(x^2-2q)z^2+(-x^2q+xq^2+3q^2)z^4+3q^3z^5-q^3z^6-4q^4z^7-2xq^4z^8}{(1-qz^2)^3(1-q^2z^3)}.$$
\end{prop}

\begin{proof}
We saw in the proof of Proposition~\ref{th:132321} that $\kra$
induces a bijection between $\S_n(132,321)$ and the set of paths
in $\D_n$ with at most two peaks. Composing it with the
complementation bijection, we get that $\pi\mapsto\kra(\bar\pi)$
is a bijection between $\S_n(123,312)$ and such set of Dyck paths.
Using Lemma~\ref{lemma:kra}, we can write
$F_{\{123,312\}}(x,q,z)=\sum x^{\tdzero(D)}q^{\tdneg(D)}z^{|D|}$,
where the sum is over Dyck paths $D$ with at most two peaks.
Again, such a $D$ can be uniquely written as $D=\u^k D'\d^k$,
where $k\ge0$ and $D'$ is either empty or a pair of adjacent
pyramids, i.e., $D'=\u^i\d^i\u^j\d^j$ with $i,j\ge1$. The idea is
to consider cases depending on the relations among $i$, $j$ and
$k$.

To enumerate Dyck paths with at most two peaks with respect to
$\tdzero(\cdot)$ and $\tdneg(\cdot)$, it is important to look at
where the tunnels of depth 0 and depth 1 occur. For convenience in
this proof, we call such tunnels \emph{frontier tunnels}, since
they determine where tunnels of depth $<0$ are: above them all
tunnels have negative depth, and below them tunnels have positive
depth. There are four possibilities according to where the
frontier tunnels of $D$ occur in the decomposition above: \ben
\item outside $D'$,
\item inside one of the pyramids of $D'$,
\item inside both pyramids of $D'$,
\item $D$ has no frontier tunnel.\een

\begin{figure}[hbt]
\epsfig{file=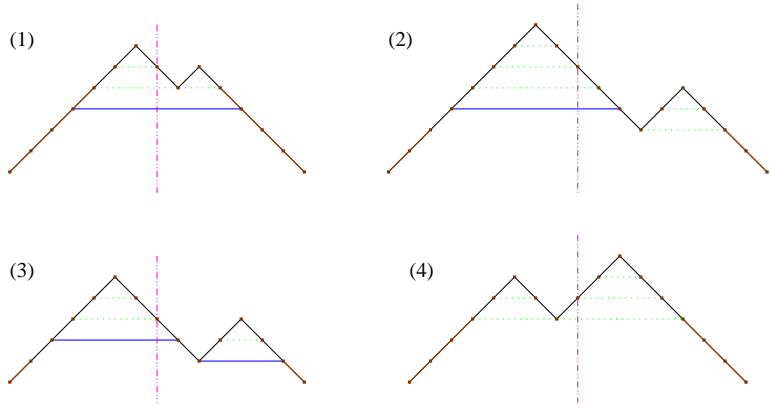,width=4in} \caption{\label{fig:casesfrontier}
Four possible locations of the frontier tunnels.}
\end{figure}

Figure~\ref{fig:casesfrontier} shows an example of each of the four cases. The frontier
tunnels (whose depth is 0 in this example) are drawn with a solid line, while
the dotted lines are the tunnels of depth $<0$.

Note that in case (4) the tunnels of negative depth are exactly
those in $D'$. We show as an example how to find the GF in case
(1). In this case, the frontier tunnel $T$ gives a decomposition
$D=A\u B\d C$ where $A=\u^m$, $C=\d^m$, $m\ge0$, and $B$ is a Dyck
path with at most two peaks, of semilength $|B|=m$ if
$\mathrm{depth}(T)=0$, and $|B|=m+1$ if $\mathrm{depth}(T)=1$. It
follows from Proposition~\ref{th:132321} that the GF for Dyck
paths with at most two peaks is $\frac{1-2z+2z^2}{(1-z)^3}$. In
the situation where $\mathrm{depth}(T)=0$, we have that
$|D|=2|B|+1$ and $\tdneg(D)=|B|$. Thus, the corresponding GF is
$$xz\cdot\frac{1-2qz^2+2q^2z^4}{(1-qz^2)^3}.$$ Similarly, in the situation where
$\mathrm{depth}(T)=1$, we have that $|D|=2|B|$ and
$\tdneg(D)=|B|$, so the corresponding GF is
$$\frac{1-2qz^2+2q^2z^4}{(1-qz^2)^3}.$$

The other cases are similar. Adding up the GFs obtained in each
case, we get the desired expression for $F_{\{123,312\}}(x,q,z)$.
\end{proof}

\begin{corollary}
$$F_{\{123,231\}}(x,q,z)
=\frac{1+xz+(x^2-2q)z^2+(-x^2q+xq+3q^2)z^4+3q^2z^5-q^3z^6-4q^3z^7-2xq^3z^8}{(1-qz^2)^3(1-qz^3)}.$$
\end{corollary}

\begin{proof}
By Lemma~\ref{lemma:sim}, we have that
$F_{\{123,231\}}(x,q,z)=F_{\{123,312\}}(\frac{x}{q},\frac{1}{q},qz)$.
\end{proof}

\subsection{h) $\{123,321\}$}

\begin{prop}\label{th:123321}
$$F_{\{123,321\}}(x,q,z)=1+xz+(x^2+q)z^2+(2xq+q^2+q)z^3+4q^2z^4.$$
\end{prop}

\begin{proof}
By a well-known result of Erd\"os and Szekeres, any permutation of
length at least 5 contains an occurrence of either 123 or 321.
This reduces the problem to counting fixed points and excedances
in permutations of length at most 4, which is trivial.
\end{proof}

\section{Triple restrictions}\label{sec:triple}

Here we consider simultaneous avoidance of any three patterns of
length~3. Using Lemma~\ref{lemma:sim}, the triplets of patterns
fall into the following equivalence classes.

\bce{\bf a) $\{123,132,213\}$ \\
b) $\{231,312,321\}$ \\
c) $\{123,132,231\}\simfe\{123,213,231\}\quad\simf\quad$ c') $\{123,132,312\}\simfe\{123,213,312\}$ \\
d) $\{132,231,321\}\simfe\{213,231,321\}\quad\simf\quad$ d') $\{132,312,321\}\simfe\{213,312,321\}$ \\
e) $\{132,213,231\}\quad\simf\quad$ e') $\{132,213,312\}$ \\
f) $\{132,231,312\}\simfe\{213,231,312\}$ \\
g) $\{123,231,312\}$ \\
h) $\{132,213,321\}$ \\
i) $\{123,132,321\}\simfe\{123,213,321\}$ \\
j) $\{123,231,321\}\quad\simf\quad$ j') $\{123,312,321\}$ \\
}\ece

It is known \cite{SS85} that the number of permutations in $S_n$
avoiding the triplets in the classes {\bf a)} and {\bf b)} is the
Fibonacci number $F_{n+1}$. The number of permutations avoiding
any of the triplets in the classes {\bf c)}, {\bf c')}, {\bf d)},
{\bf d')}, {\bf e)}, {\bf e')}, {\bf f)}, {\bf g)} and {\bf h)} is
$n$. The cases of the triplets {\bf i)}, {\bf j)} and {\bf j')}
are trivial, because they are avoided only by permutations of
length at most 4.

In terms of generating functions, when we substitute $x=q=1$ in
$F_\Sigma(x,q,z)$ where $\Sigma$ is a triplet from one of the
classes between {\bf a)} and {\bf g)}, we get
$$F_\Sigma(1,1,z)=\sum_{n\ge0}F_{n+1}z^n=\frac{1}{1-z-z^2}.$$
If $\Sigma$ is any triplet from the classes between {\bf c)} and
{\bf h)}, we get
$$F_\Sigma(1,1,z)=\sum_{n\ge0}n z^n=\frac{1-z+z^2}{(1-z)^2}.$$

The following theorem gives all the generating functions of
permutations avoiding any triplet of patterns of length~3.

\begin{theorem}\label{th:triple}

{\bf a)}
$$F_{\{123,132,213\}}(x,q,z)=\frac{1+xz+(x^2-q)z^2+(-xq+q^2+q)z^3-x^2qz^4}{(1+qz^2)(1-3qz^2+q^2z^4)}$$

{\bf b)} $$F_{\{231,312,321\}}(x,q,z)=\frac{1}{1-xz-qz^2}$$

{\bf c)}
$$F_{\{123,132,231\}}(x,q,z)=F_{\{123,213,231\}}(x,q,z)=\frac{1+xz+(x^2-q)z^2+qz^3+(-x^2q+xq+q^2)z^4}{(1-qz^2)^2}$$

{\bf c')}
$$F_{\{123,132,312\}}(x,q,z)=F_{\{123,213,312\}}(x,q,z)=\frac{1+xz+(x^2-q)z^2+q^2z^3+(-x^2q+xq^2+q^2)z^4}{(1-qz^2)^2}$$

{\bf d)}
$$F_{\{132,231,321\}}(x,q,z)=F_{\{213,231,321\}}(x,q,z)=\frac{1-z+qz^2}{(1-z)(1-xz)}$$

{\bf d')}
$$F_{\{132,312,321\}}(x,q,z)=F_{\{213,312,321\}}(x,q,z)=\frac{1-qz+qz^2}{(1-xz)(1-qz)}$$

{\bf e)}
$$F_{\{132,213,231\}}(x,q,z)=\frac{1-z-qz^2+2qz^3+(-x^2q+q^2-xq)z^4+(x^2q-2q^2)z^5+xq^2z^6}{(1-z)(1-xz)(1-qz^2)^2}$$

{\bf e')}
$$F_{\{132,213,312\}}(x,q,z)=\frac{1-qz-qz^2+2q^2z^3+(-x^2q-xq^2+q^2)z^4+(x^2q^2-2q^3)z^5+xq^3z^6}{(1-xz)(1-qz)(1-qz^2)^2}$$

{\bf f)}
$$F_{\{132,231,312\}}(x,q,z)=F_{\{213,231,312\}}(x,q,z)=\frac{1+xqz^3}{(1-xz)(1-qz^2)}$$

{\bf g)}
$$F_{\{123,231,312\}}(x,q,z)=\frac{1+xz+(x^2-q)z^2+xqz^3+q^2z^4}{(1-qz^2)^2}$$

{\bf h)}
$$F_{\{132,213,321\}}(x,q,z)=\frac{1-(1+q)z+2qz^2-xqz^3}{(1-z)(1-xz)(1-qz)}$$

{\bf i)}
$$F_{\{123,132,321\}}(x,q,z)=F_{\{123,213,321\}}(x,q,z)=1+xz+(x^2+q)z^2+(xq+q^2+q)z^3+q^2z^4$$

{\bf j)}
$$F_{\{123,231,321\}}(x,q,z)=1+xz+(x^2+q)z^2+(2xq+q)z^3+q^2z^4$$

{\bf j')}
$$F_{\{123,312,321\}}(x,q,z)=1+xz+(x^2+q)z^2+(2xq+q^2)z^3+q^2z^4$$
\end{theorem}

\begin{proof}
Throughout this proof we will use the bijection
$\kra:\S_n(132)\longrightarrow\D_n$ described in
Subsection~\ref{sec:bij}.

{\bf a)} As in the proof of Proposition~\ref{th:123132}, we have
that $\pi\in\S_n(132)$ avoids 123 if and only if the Dyck path
$\kra(\pi)$ has height at most 2. Similarly, from the proof of
Proposition~\ref{th:132213}, $\pi$ avoids 213 if and only if
$\kra(\pi)$ is a sequence of pyramids. Thus, $\kra$ induces a
bijection between $\S_n(123,132,213)$ and
$\Pyr^{\le2}_n:=\Pyr^{\le2}\cap\D_n$, where $\Pyr^{\le2}$ denotes
the set of sequences of pyramids of height at most 2. By
Lemma~\ref{lemma:kra},
$$F_{\{123,132,213\}}(x,q,z)=\sum_{D\in\Pyr^{\le2}}
x^{\ct(D)}q^{\rt(D)}z^{|D|}.$$

To count centered and right tunnels, we distinguish cases
according to which steps are the middle steps of $D$. A path in
$\Pyr^{\le2}$ of height~0 at the middle can be split in two
elements of $\Pyr^{\le2}$ of equal length, only the right one
producing right tunnels. Since the number of $D\in\Pyr^{\le2}_n$
is $F_{n+1}$, the GF for paths of height~0 at the middle is
$$\sum_{n\ge0}F_{m+1}^2 q^m
z^{2m}=\frac{1-qz^2}{(1+qz^2)(1-3qz^2+q^2z^4)}.$$ Multiplying this
expression by $xz$ (resp. by $x^2z^2$) we obtain the GF for paths
in $\Pyr^{\le2}$ having in the middle a centered pyramid of
height~1 (resp. of height 2).

\begin{figure}[hbt]
\epsfig{file=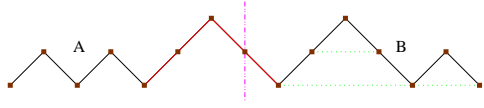,width=2.5in} \caption{\label{fig:pyr2} A
sequence of pyramids of height at most 2.}
\end{figure}

Paths $D\in\Pyr^{\le2}$ whose two middle steps are $\d\d$ can be
written as $D=A\u\u\d\d B$, where $A,B\in\Pyr^{\le2}$ and
$|B|=|A|+1$ (see Figure~\ref{fig:pyr2}). Thus, the corresponding GF is
$$\sum_{n\ge1}F_m F_{m+1} q^m z^{2m+1}=\frac{qz^3}{(1+qz^2)(1-3qz^2+q^2z^4)}.$$
By symmetry, multiplying this expression by $q$ we get the GF for
paths whose two middle steps are $\u\u$.

Adding up all the cases, we get the desired GF
$$F_{\{123,132,213\}}(x,q,z)=\frac{(1+xz+x^2z^2)(1-qz^2)}{(1+qz^2)(1-3qz^2+q^2z^4)}+\frac{(q+q^2)z^3}{(1+qz^2)(1-3qz^2+q^2z^4)}.$$

\ms

{\bf b)} Using the same reasoning as in {\bf a)}, we have that
$\pi\mapsto\kra(\bar\pi)$ induces a bijection between
$\S_n(231,312,321)$ and $\Pyr^{\le2}_n$. Now,
Lemma~\ref{lemma:kra} implies that
$$F_{\{231,312,321\}}(x,q,z)=\sum_{D\in\Pyr^{\le2}}
x^{\tdzero(D)}q^{\tdneg(D)}z^{|D|}.$$ Each pyramid of height~1
produces a tunnel of depth 0, and each pyramid of height 2 creates
a tunnel of depth $<0$. Therefore,
$$F_{\{231,312,321\}}(x,q,z)=\frac{1}{1-xz-qz^2}.$$

\ms

{\bf c)} We saw in the proof of Proposition~\ref{th:132231} that
$\pi\in\S_n(132)$ avoids 231 if and only if the Dyck path
$\kra(\pi)$ does not have any two consecutive up-steps after the
first down-step. Therefore, $\kra$ induces a bijection between
$\S_n(123,132,312)$, and paths in $\D_n$ with the above condition
and height at most 2. Such paths (except the empty one) can be
expressed uniquely as $D=\u A\d B$, where $A$ and $B$ are
sequences of hills (i.e, they have the form $(\u\d)^k$ for some
$k\ge0$). Lemma~\ref{lemma:kra} reduces the problem to enumerating
centered tunnels and right tunnels on these paths.

If $B$ is empty, $D=\u A\d$ has a centered tunnel at height~0. The
contribution of paths of this kind to our GF is
$\frac{xz}{1-qz^2}$ for $|A|$ even, and $\frac{x^2z^2}{1-qz^2}$
for $|A|$ odd.

Assume now that $|A|<|B|$, so that $A$ is within the left half of
$D=\u A\d B$. If the middle of $D$ is at height~0, then $D$ is
determined by the length of $A$ and the number of hills in $B$ to
the left of the middle. Thus, the contribution of this subset to
the GF is
$$\frac{qz^2}{(1-qz^2)^2}.$$ Multiplying this expression by $xz$
gives the GF for paths whose midpoint is on top of a hill of $B$.

\begin{figure}[hbt]
\epsfig{file=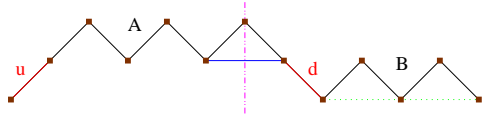,width=2.5in} \caption{\label{fig:uadb} An
example with $|A|=3$ and $|B|=2$.}
\end{figure}

It remains the case in which $|A|\ge|B|>0$. If $|A|-|B|$ is even,
the contribution of these paths to the GF is
$$z\cdot\frac{qz^2}{1-qz^2}\cdot\frac{1}{1-qz^2},$$
where the last factor counts how larger $A$ is than $B$. If
$|A|-|B|$ is odd, the corresponding GF is
$$z\cdot\frac{qz^2}{1-qz^2}\cdot\frac{xz}{1-qz^2},$$
since in this case there is a centered tunnel of height~1 inside
$A$ (see Figure~\ref{fig:uadb}).

Summing up all the cases, we get
$$F_{\{123,132,231\}}(x,q,z)=1+\frac{xz+x^2z^2}{1-qz^2}+
\frac{(1+xz)qz^2}{(1-qz^2)^2}+\frac{qz^3(1+xz)}{(1-qz^2)^2}.$$

\ms

{\bf c')} By Lemma~\ref{lemma:sim}, we have that
$F_{\{123,132,312\}}(x,q,z)=F_{\{123,132,231\}}(\frac{x}{q},\frac{1}{q},qz)$,
so the formula follows from part {\bf c)}.

\ms

{\bf d)} As in the proof of Proposition~\ref{th:132231}, we use
that $\pi\in\S_n(132)$ avoids 231 if and only if the Dyck path
$\kra(\pi)$ does not have any two consecutive up-steps after the
first down-step. Besides, as in Proposition~\ref{th:132321},
$\pi\in\S_n(132)$ avoids 321 if and only if $\kra(\pi)$ has at
most two peaks. Thus, $\pi\in\S_n(132,231,321)$ if and only if
$\kra(\pi)\in\D_n$ has the form $\u^k B\d^k$, where $B$ is either
empty or is a pair of pyramids, the second of height~1. Fixed
points and excedances of $\pi$ are mapped to centered tunnels and
right tunnels of $\kra(\pi)$ respectively, by
Lemma~\ref{lemma:kra}. Thus, $F_{\{123,132,312\}}(x,q,z)$ equals
the GF enumerating centered and right tunnels in these paths.

If $B$ is not empty, the contribution of the first pyramid is
$\frac{z}{1-z}$, and the second pyramid contributes $qz$. Centered
tunnels come from the steps outside $B$. Hence,
$$F_{\{123,132,312\}}(x,q,z)=\frac{1}{1-xz}\left(1+\frac{z}{1-z}\cdot qz\right).$$

\ms

{\bf d')} It follows from part {\bf d)} and Lemma~\ref{lemma:sim}.

\ms

{\bf e)} Let $\pi\in\S_n(132)$. We have seen that the condition
that $\pi$ avoids 213 translates into $\kra(\pi)$ being a sequence
of pyramids. The additional restriction of $\pi$ avoiding 231
implies that all but the first pyramid of the sequence $\kra(\pi)$
must have height~1. Thus, by Lemma~\ref{lemma:kra},
$F_{\{132,213,231\}}(x,q,z)$ can be obtained enumerating centered
and right tunnels in paths of the form $D=AB$, where $A$ is any
pyramid and $B$ is a sequence of hills.

The contribution of such paths when $B$ is empty is just
$\frac{1}{1-xz}$. Assume now that $B$ is not empty. If $|A|>|B|$,
the corresponding contribution is
$$\frac{qz^2}{1-qz^2}\cdot\frac{z}{1-z},$$ where the second factor
counts how larger $A$ is than $B$. It remains the case
$|A|\le|B|$, in which $A$ is within the left half of $D$. If the
middle of $D$ is at height~0, then $D$ is determined by the length
of $A$ and the number of hills in $B$ to the left of the middle.
Thus, the contribution of this subset to the GF is
$$\frac{qz^2}{(1-qz^2)^2}.$$ Multiplying this expression by $xz$
gives the GF for paths whose midpoint is on top of a hill of $B$.

\begin{figure}[hbt]
\epsfig{file=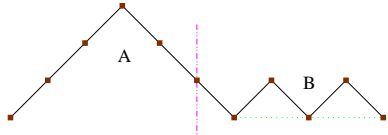,width=2in} \caption{\label{fig:pyrhill} A
pyramid followed by a sequence of hills.}
\end{figure}

Summing all this up, we get
$$F_{\{132,213,231\}}(x,q,z)=\frac{1}{1-xz}+\frac{qz^3}{(1-z)(1-qz^2)}
+\frac{(1+xz)qz^2}{(1-qz^2)^2}.$$

\ms

{\bf e')} It follows from part {\bf e)} and Lemma~\ref{lemma:sim}.

\ms

{\bf f)} Reasoning as in the proof of {\bf e)}, we see that
$\pi\mapsto\kra(\bar\pi)$ induces a bijection between
$\S_n(123,231,312)$ and the subset of paths in $D_n$ consisting of a
pyramid followed by a sequence of hills. By Lemma~\ref{lemma:kra},
it is enough to enumerate these paths according to the statistics $\tdzero(\cdot)$ and
$\tdneg(\cdot)$. If the path is not empty, the first pyramid contributes
$\frac{xz}{1-qz^2}$ if it has odd size (since then it contains a
tunnel of depth 0) and $\frac{qz^2}{1-qz^2}$ if it has even size.
The sequence of hills contributes $\frac{1}{1-xz}$. Therefore,
$$F_{\{132,231,312\}}(x,q,z)=1+\frac{xz+qz^2}{1-qz^2}\cdot\frac{1}{1-xz}.$$

\ms

{\bf g)} Let $\pi\in\S_n(132)$. We have seen that $\pi$ avoids 213
if and only if $\kra(\pi)$ is a sequence of pyramids, and that
$\pi$ avoids 321 if and only if $\kra(\pi)$ has at most two peaks.
In other words, $\kra$ induces a bijection between
$\S_n(132,213,321)$ and the subset of paths in $\D_n$ that are a
sequence of at most two pyramids. Composing with the
complementation operation, we have that $\pi\in\S_n(123,231,312)$
if and only if $\kra(\bar\pi)$ is in that subset. Now,
Lemma~\ref{lemma:kra} implies that $F_{\{123,231,312\}}$ can be
obtained enumerating sequences of at most 2 pyramids according to
$\tdzero(\cdot)$ and $\tdneg(\cdot)$. Each pyramid contributes
$\frac{xz}{1-qz^2}$ if it has odd size and $\frac{qz^2}{1-qz^2}$
if it has even size. Thus,
$$F_{\{123,231,312\}}(x,q,z)=1+\frac{xz+qz^2}{1-qz^2}+\left(\frac{xz+qz^2}{1-qz^2}\right)^2.$$

\ms

{\bf h)} We have shown in the proof of {\bf g)} that
$\pi\in\S_n(132,213,321)$ if and only if $\kra(\pi)$ is a sequence
of at most two pyramids. Using Lemma~\ref{lemma:kra}, it is enough
to enumerate centered tunnels and right tunnels in such paths. The
contribution of paths with exactly two pyramids is
$$\frac{z}{1-z}\cdot\frac{qz}{1-qz},$$
since only the one on the right gives right tunnels. Centered
tunnels appear when there is only one pyramid. Thus we obtain
$$F_{\{132,213,321\}}(x,q,z)=\frac{1}{1-xz}+\frac{qz^2}{(1-z)(1-qz)}.$$

\ms

{\bf i, j, j')} These cases are trivial because only permutations
of length at most 4 can avoid 123 and 321 simultaneously.
\end{proof}

\ms

After having studied all the cases of double and triple
restrictions, the next step is to consider restrictions of higher
multiplicity. However, for $\Sigma\subseteq\S_3$, $|\Sigma|\ge4$,
the sets $S_n(\Sigma)$ are very easy to describe (see for example
\cite{SS85}), and the distribution of fixed points and excedances
is trivial. In particular, in these cases
$|S_n(\Sigma)|\in\{0,1,2\}$ for all $n$.

\section{Pattern-avoiding involutions}\label{sec:invol}

Recall that a permutation $\pi$ is an involution if
$\pi=\pi^{-1}$. In terms of the array representation of $\pi$,
this condition is equivalent to the array being symmetric with
respect to the main diagonal. Denote by $\I_n$ the set of
involutions of length $n$. In this section we consider the
distribution of the statistics $\fp(\cdot)$ and $\exc(\cdot)$ in
involutions avoiding any subset of patterns of length~3.

For any $\pi\in\S_n$, it is clear that
$\fp(\pi)+\exc(\pi)+\exc(\pi^{-1})=n$ (each cross in the array of
$\pi$ is either on, to the right of, or to the left of the main
diagonal). Thus, if $\pi\in\I_n$, then
$\exc(\pi)=\frac{1}{2}(n-\fp(\pi))$, so the number of excedances
is determined by the number of fixed points. Therefore, it is
enough here to consider only the statistic `number of fixed
points' in pattern-avoiding involutions.

For any set of patterns $\Sigma$, let
$\I_n(\Sigma):=\I_n\cap\S_n(\Sigma)$, and let
$i_n^k(\Sigma):=|\{\pi\in\I_n(\Sigma):\fp(\pi)=k\}|$. Define
$$G_\Sigma(x,z):=\sum_{n\ge0}\sum_{\pi\in\I_n(\Sigma)}x^{\fp(\pi)}z^n.$$
By the reasoning above,
$\sum_{n\ge0}\sum_{\pi\in\I_n(\Sigma)}x^{\fp(\pi)}q^{\exc(\pi)}z^n=G_\Sigma(xq^{-1/2},zq^{1/2}$).

Clearly, $\wh\pi$ is an involution if and only if $\pi$ is an
involution. Therefore, from Lemma~\ref{lemma:refl} we get the
following.

\begin{lemma}
\label{lemma:siminv} Let $\Sigma$ be any set of permutations. We
have \ben
\item $G_{\wh\Sigma}(x,z)=G_{\Sigma}(x,z)$,
\item
$G_{\Sigma^{-1}}(x,z)=G_{\Sigma}(x,z)$. \een
\end{lemma}

The main reason our techniques for studying statistics in
pattern-avoiding permutations can be applied to the case of
involutions is stated in the following lemma.

\begin{lemma} \label{lemma:krasym}
Let $\pi\in\S_n(132)$ and let $D=\kra(\pi)\in\D_n$. Then,
$$\pi\mbox{ is an involution}\Longleftrightarrow\kra(\pi)\mbox{ is symmetric}.$$
\end{lemma}
\begin{proof}
The array of crosses representing $\pi^{-1}$ is obtained from the
one of $\pi$ by reflection over the main diagonal. Therefore, from
the description of the bijection $\kra$ given in
Subsection~\ref{sec:bij}, we have that $\kra(\pi^{-1})=D^\ast$. It
follows that $\pi$ is an involution if and only if $D=D^\ast$,
which is equivalent to $D$ being a symmetric Dyck path.
\end{proof}

\subsection{Single restrictions}

It is known \cite{SS85} that for $\sigma\in\{123,132,213,321\}$,
$|\I_n(\sigma)|=\binom{n}{\lfloor n/2\rfloor}$, and that for
$\sigma\in\{231,312\}$, $|\I_n(\sigma)|=2^{n-1}$. From
Lemma~\ref{lemma:siminv} it follows that for all $k\ge0$,
$i_n^k(132)=i_n^k(213)$ and $i_n^k(231)=i_n^k(312)$. It is shown
in~\cite{DRS} (see also~\cite{EliPak} for a bijective proof) a
that in fact $i_n^k(132)=i_n^k(321)$. So, for single restrictions
there are three cases to consider.

\begin{theorem}[\cite{GM02,DRS}] Let $n\ge1$, $k\ge0$. We have
\ben \item $i_n^0(123)=i_n^2(123)=\left\{\bt{ll} $\binom{n-1}{\frac{n}{2}}$ & if $n$ is even, \\
$0$ & if $n$ is odd, \et\right.$ \\
$i_n^1(123)=\left\{\bt{ll} $\binom{n}{\frac{n-1}{2}}$ & if $n$ is odd, \\
$0$ & if $n$ is even, \et\right.$ \\  $i_n^k(123)=0$ if $k\ge3$.
\item $i_n^k(132)=i_n^k(213)=i_n^k(321)=\left\{\bt{ll} $\frac{k+1}{n+1}\binom{n+1}{\frac{n-k}{2}}$ & if $n-k$ is even, \\
$0$ & if $n-k$ is odd. \et\right.$
\item $G_{231}(x,z)=G_{312}(x,z)=\dfrac{1-z^2}{1-xz-2z^2}$.
\een
\end{theorem}

\begin{proof}
(1) Clearly a 123-avoiding permutation cannot have more than two
fixed points. On the other hand, if $\pi\in\I_n$, we have
$\fp(\pi)=n-2\exc(\pi)$, which explains that $i_n^k(123)=0$ if
$n-k$ is odd. This implies that for odd $n$, $\fp(\pi)=1$ for all
$\pi\in\I_n$, so
$i_n^1(123)=|\I_n(123)|=\binom{n}{\frac{n-1}{2}}$. For even $n$,
all we have to show is that $i_n^0(123)=i_n^2(123)$.

The bijection $\brs:\S_n(123)\longrightarrow\D_n$ described in
Subsection~\ref{s:123} has the property that $\pi\in\I_n(123)$ if
and only if $\brs(\pi)$ is a symmetric Dyck path. If $n$ is even,
involutions $\pi\in\I_n$ with $\fp(\pi)=2$ are mapped to symmetric
Dyck paths with a peak in the middle, and those with $\fp(\pi)=0$
are mapped to symmetric Dyck paths with a valley in the middle. We
can establish a bijection between these two sets of Dyck paths
just by changing the middle peak $\u\d$ into a middle valley
$\d\u$ (this can always be done because the height at the middle
of a Dyck path of even semilength is always even, so it cannot be
1). This proves that $i_n^0(123)=i_n^2(123)$, and in particular it
equals $\frac{1}{2}|\I_n(123)|=\binom{n-1}{\frac{n}{2}}$.

(2) We use the bijection $\kra:\S_n(132)\longrightarrow\D_n$,
which by Lemma~\ref{lemma:krasym} restricts to a bijection between
$\I_n(132)$ and $\Ds$. Thus, by Lemma~\ref{lemma:kra},
$G_{123}(x,z)$ can be expressed as $\sum_{D\in\Ds}
x^{\ct(D)}z^{|D|}$, where the sum is over all symmetric Dyck
paths. But the number of centered tunnels of a symmetric Dyck path
is just its height at the middle. Therefore, taking only the first
half of the path, $i_n^k(123)$ counts the number of paths from
$(0,0)$ to $(n,k)$ never going below the $x$-axis, which equals
the ballot number given in the theorem.

(3) Consider the bijection \bt{ccc} $\S_n(312)$ &
$\longleftrightarrow$ & $\D_n$ \\ $\pi$ & $\mapsto$ &
$\kra(\bar\pi)$ \et. Then $\pi\in\I_n(312)$ if and only if
$\kra(\bar\pi)$ is a sequence of pyramids. In fact, it turns out
\cite{SS85} that $\I_n(312)=\I_n(231)=\S_n(231,312)$. By
Lemma~\ref{lemma:kra}, fixed points of $\pi$ are mapped to tunnels
of depth 0 of $\kra(\bar\pi)$, which are produced by pyramids of
odd size. Thus, as in Proposition~\ref{th:231312},
$$G_{312}(x,z)=\frac{1}{1-\frac{xz+z^2}{1-z^2}}.$$
\end{proof}

\subsection{Multiple restrictions}
\begin{theorem}
{\bf a)}
$$G_{\{123,132\}}(x,z)=G_{\{123,213\}}(x,z)=\frac{1+xz+(x^2-1)z^2}{1-2z^2}$$

{\bf b)}
$$G_{\{231,321\}}(x,z)=G_{\{312,321\}}(x,z)=\frac{1}{1-xz-z^2}$$

{\bf c)} $$G_{\{132,213\}}(x,z)=\frac{1-z^2}{(1-xz)(1-2z^2)}$$

{\bf d)} $$G_{\{231,312\}}(x,z)=\frac{1-z^2}{1-xz-2z^2}$$

{\bf e)}
$$G_{\{132,231\}}(x,z)=G_{\{213,231\}}(x,z)=G_{\{132,312\}}(x,z)=G_{\{213,312\}}(x,z)=\frac{1+xz^3}{(1-xz)(1-z^2)}$$

{\bf f)}
$$G_{\{132,321\}}(x,z)=G_{\{213,321\}}(x,z)=\frac{1}{(1-xz)(1-z^2)}$$

{\bf g)}
$$G_{\{123,231\}}(x,z)=G_{\{123,312\}}(x,z)=\frac{1+xz+(x^2-1)z^2+xz^3+z^4}{(1-z^2)^2}$$

{\bf h)} $$G_{\{123,321\}}(x,z)=1+xz+(x^2+1)z^2+2xz^3+2z^4$$
\end{theorem}

\begin{proof}
All the equalities between different $G_\Sigma$ follow trivially
from Lemma~\ref{lemma:siminv}. To find expressions for these GFs,
the idea is to use again the same bijections as in
Section~\ref{sec:double}, between permutations avoiding two
patterns of length~3 and certain subclasses of Dyck paths. The
main difference is that here we will have to deal only with
symmetric Dyck paths, as a consequence of
Lemma~\ref{lemma:krasym}.

{\bf a)} From the proof of Proposition~\ref{th:123132} and
Lemma~\ref{lemma:krasym}, we have that $\kra$ restricts to a
bijection between $\I_n(123,132)$ and symmetric Dyck paths
$D\in\D_n$ of height at most 2. By Lemma~\ref{lemma:kra}, $\kra$
maps fixed points to centered tunnels, so all we have to do is
count elements $D\in\Ds$ of height at most 2 according to the
number of centered tunnels. Such a $D$ can be uniquely written as
$D=ABC$, where $A=C^\ast\in\D^{\le2}$ and $B$ is either empty or
has the form $B=\u B_1\d$, where $B_1$ is a sequence of hills. If
$|B_1|$ is even (resp. odd), then $D$ has one (resp. two) centered
tunnels, so the contribution of $B$ is $1+\frac{(1+xz)xz}{1-z^2}$.
The contribution of $A$ and $C$ is $\frac{1-z^2}{1-2z^2}$. The
product of these two quantities gives the expression for
$G_{\{123,132\}}(x,z)$.

{\bf b)} We have seen (\cite{SS85}) that
$\I_n(231)=\S_n(231,312)$. Therefore,
$\I_n(231,321)=\S_n(231,312,321)$. This case was treated in
Theorem~\ref{th:triple} {\bf b)}.

{\bf c)} From the proof of Proposition~\ref{th:132213} and
Lemma~\ref{lemma:krasym}, we have that $\kra$ gives a bijection
between $\I_n(132,213)$ and symmetric sequences of pyramids
$D\in\Pyr_n$, and that it maps fixed points of the permutation to
centered tunnels of the Dyck path. Such a $D$ can be written
uniquely as $D=ABC$, where $A=C^\ast\in\Pyr$, and $B$ is either
empty or a pyramid. The contribution of $B$ is $\frac{1}{1-xz}$,
whereas $A$ and $C$ contribute $\frac{1-z^2}{1-2z^2}$. Multiplying
these two expressions we get a formula for $G_{\{132,213\}}(x,z)$.

{\bf d)} Again, $\I_n(231)=\S_n(231,312)$ implies that
$\I_n(231,312)=\S_n(231,312)$, which has been considered in
Proposition~\ref{th:231312}.

{\bf e)} We have that $\I_n(132,231)=\S_n(132,231,312)$, so the
fromula follows from Theorem~\ref{th:triple} {\bf f)}.

{\bf f)} From the proof of Proposition~\ref{th:132321} and
Lemma~\ref{lemma:krasym}, we have that $\kra$ gives a bijection
between $\I_n(132,321)$ and symmetric paths $D\in\D_n$ with at
most two peaks. Counting centered tunnels in such paths is very
easy, since they have the form $D=\u^k B\d^k$, where $k\ge0$ and
$B$ is either empty or a pair of identical pyramids. The
contribution of $B$ is $\frac{1}{1-z^2}$, whereas the rest
contributes $\frac{1}{1-xz}$.

{\bf g)} We have that $\I_n(123,231)=\S_n(123,231,312)$, so the
formula follows from Theorem~\ref{th:triple} {\bf g)}.

{\bf h)} It is trivial since $\S_n(123,321)=\emptyset$ for
$n\ge5$.
\end{proof}

The case of involutions avoiding simultaneously three or more
patterns of length~3 is very easy and does not involve any new
idea, so we omit it here.

\section{Final remarks}

Looking at the results of this paper, one can observe that the GFs
$F_\Sigma(x,q,z)$ that we have obtained for $\Sigma\subseteq\S_3$
are all rational functions when $|\Sigma|\ge2$. This contrasts with
the fact that they are
not rational when $|\Sigma|=1$, since in that case
$F_\Sigma(1,1,z)=\frac{1-\sqrt{1-4z}}{2z}=\C(z)$. For the case of
involutions, all the GFs $G_\Sigma(x,z)$ for $\Sigma\subseteq\S_3$
are rational except when
$\Sigma\in\{\{123\},\{132\},\{213\},\{321\}\}$.

\ms

Regarding possible extensions of this work, it would be
interesting to find a generating function for fixed points and
excedances in 123-avoiding permutations, the only case of patterns
of length~3 that remains unsolved. Even for the enumeration of
fixed points in these permutations, we expect that an expression
simpler than the one in Theorem~\ref{th:messy123} can be given.

Another further direction of research would consist in describing
the cycle structure of pattern-avoiding permutations. Using the
same bijective techniques from this paper, one can easily derive
generating functions for the \emph{augmented cycle index} of
permutations in $\S_n(231,312)$, in $\S_n(231,321)$ and in
$\S_n(132,321)$. However, it is not clear whether for permutations
avoiding other subsets of patterns of length~3, the distribution
of the cycle type has a simple description.

\end{document}